\documentclass[preprint, 3p, times, 12pt]{elsarticle}

\usepackage{amssymb}
\usepackage{amsmath}
\usepackage{graphicx}
\usepackage{amsmath}
\usepackage{amssymb}
\usepackage{amsthm}
\usepackage[colorlinks=true]{hyperref}
\usepackage[linesnumbered,ruled,vlined]{algorithm2e}
\usepackage{caption}
\usepackage{subcaption}
\usepackage[capitalise]{cleveref}
\usepackage{multirow}
\usepackage{booktabs}

\newtheorem*{remark}{Remark}

\SetAlgoHangIndent{0pt}
\DeclareMathOperator*{\argmin}{arg\,min}

\journal{Journal of Computational Physics}

\begin{document}

\begin{frontmatter}

\title{PACMANN: Point Adaptive Collocation Method for Artificial Neural Networks}

\author[a]{Coen Visser}

\author[b]{Alexander Heinlein\corref{corresponding_author}}
\ead{A.Heinlein@tudelft.nl}

\author[a,c]{Bianca Giovanardi}

\affiliation[a]{organization={Faculty of Aerospace Engineering},
            addressline={Delft University of Technology},
            % city={},
            state={2629 HS Delft},
            country={The Netherlands}}

\affiliation[b]{organization={Delft Institute of Applied Mathematics, Faculty of Electrical Engineering, Mathematics and Computer Science},
            addressline={Delft University of Technology},
            % city={},
            state={2628 CD Delft},
            country={The Netherlands}}

\affiliation[c]{organization={Delft Institute for Computational Science and Engineering},
            addressline={Delft University of Technology},
            % city={},
            state={2628 CD Delft},
            country={The Netherlands}}

\cortext[corresponding_author]{Corresponding author}

\begin{abstract}
Physics-Informed Neural Networks (PINNs) have emerged as a tool for approximating the solution of Partial Differential Equations (PDEs) in both forward and inverse problems. PINNs minimize a loss function which includes the PDE residual determined for a set of collocation points. Previous work has shown that the number and distribution of these collocation points have a significant influence on the accuracy of the PINN solution. Therefore, the effective placement of these collocation points is an active area of research. Specifically, available adaptive collocation point sampling methods have been reported to scale poorly in terms of computational cost when applied to high-dimensional problems. In this work, we address this issue and present the Point Adaptive Collocation Method for Artificial Neural Networks (PACMANN). PACMANN incrementally moves collocation points toward regions of higher residuals using gradient-based optimization algorithms guided by the gradient of the PINN loss function, that is, the squared PDE residual. We apply PACMANN for forward and inverse problems, and demonstrate that this method matches the performance of state-of-the-art methods in terms of the accuracy/efficiency tradeoff for the low-dimensional problems, while outperforming available approaches for high-dimensional problems. Key features of the method include its low computational cost and simplicity of integration into existing physics-informed neural network pipelines. The code is available at \url{https://github.com/CoenVisser/PACMANN}.
\end{abstract}

% \begin{highlights}
% \end{highlights}


\end{frontmatter}

\section{Introduction}
\label{sec:Introduction}
Physics-Informed Neural Networks (PINNs) build upon the ability of deep neural networks to serve as universal function approximators, as established by Cybenko~\cite{Cybenko1989ApproximationFunction} and Hornik et al.~\cite{Hornik1989MultilayerApproximators} in 1989. Based on these findings, several methods were developed to solve Ordinary Differential Equations (ODEs) and Partial Differential Equations (PDEs) using neural networks, originally proposed by  \cite{Dissanayake1994Neural-network-basedEquations,Lagaris1998ArtificialEquations}. Supported by these developments and recent advances in computational tools, notably automatic differentiation in 2015 \cite{Baydin2015AutomaticSurvey}, Raissi et al.~\cite{Raissi2019Physics-informedEquations} proposed the name and framework of \emph{Physics-Informed Neural Networks} and their use to approximate the solution of PDEs in both forward and inverse problems; their work was published in 2019. Since then, PINNs have been applied in a variety of fields \cite{Cuomo2022ScientificNext, Karniadakis2021Physics-informedLearning}, such as fluid dynamics \cite{Mao2020Physics-informedFlows, Cai2021FlowNetworks, Jin2021NSFnetsEquations}, heat transfer \cite{Cai2021Physics-informedProblems, AminiNiaki2021Physics-informedManufacture}, material sciences \cite{Shukla2020Physics-informedCracks, Zhang2022AnalysesNetworks}, and electromagnetism \cite{Kovacs2022ConditionalNetworks, Son2023AMotor}.

In their classic form, PINNs approximate the solution of differential equations by minimizing a loss function incorporating boundary conditions, initial conditions, and the PDE residual sampled over a set of collocation points. In 2020, Mao et al.~\cite{Mao2020Physics-informedFlows} explored the impact of collocation point placement on prediction accuracy for solutions exhibiting discontinuities. They demonstrated that, when discontinuities are known a priori, manually increasing the density of points near these regions improves prediction accuracy.
To address scenarios where solution features are unknown before training, adaptive algorithms for collocation point selection were developed. For instance, Lu et al.~\cite{Lu2021DeepXDE:Equations} introduced Residual-based Adaptive Refinement (RAR) in 2021, the first adaptive sampling algorithm, which places additional points in regions with the largest PDE residuals. RAR improves prediction accuracy, capturing features such as a discontinuity better than a static grid for the Burgers' equation. Instead of sampling points only where the residual is the largest, Nabian et al.~\cite{Nabian2021EfficientSampling} proposed to randomly resample all points in the domain based on a Probability Density Function (PDF) proportional to the loss function. This approach samples a higher density of points in high residual areas, resulting in accelerated convergence of the PINN. Based upon these studies, Wu et al.~\cite{Wu2023ANetworks} presented two additional sampling algorithms, Residual-based Adaptive Distribution (RAD) and Residual-based Adaptive Refinement with Distribution (RAR-D). In RAD, all collocation points are resampled using a PDF defined by a nonlinear function of the PDE residual. RAR-D is a combination of RAR and RAD, where collocation points are sampled in addition to the existing ones according to the same probability density function used by RAD. Both approaches lead to a higher prediction accuracy, specifically for PDEs with complex solutions due to, for example, steep gradients. Moreover, RAD was found to outperform the method proposed by Nabian et al.~\cite{Nabian2021EfficientSampling}.

Several adaptations to the aforementioned PDF-based sampling algorithms have been proposed. Guo et al.~\cite{Guo2024TCAS-PINN:Method}, for example, propose an adaptive causal sampling method, which decomposes the domain into subdomains where the ratio of points sampled in each subdomain is based on the PDE residual and a temporal weight, ensuring temporal causality. This approach was found to enhance the prediction accuracy and computational efficiency of PINNs in problems with nonlinear PDEs containing higher-order derivatives. Furthermore, Mao et al.~\cite{Mao2023Physics-informedSolutions} consider the gradient of the solution by sampling additional points in subdomains with large residuals and large solution gradients. Liu et al.~\cite{Liu2024AnPINNs} propose to add points with large residual gradients to the set of collocation points used for training. Both Mao et al. and Liu et al. report an improvement in accuracy for problems with solutions exhibiting steep gradients.

While the aforementioned collocation point sampling methods have proven effective in low-dimensional problems, these approaches to resampling are computationally expensive for high-dimensional problems, as reported by Wu et al.~\cite{Wu2023ANetworks}.
Specifically, RAD or RAR are computationally expensive for these problems due to the cost of evaluating the residual at a sufficiently large number of points, either to construct the probability density function or to identify additional points for inclusion in the training process. Other approaches have been proposed to sample collocation points in high-dimensional problems. For instance, Tang et al.~\cite{Tang2021DAS-PINNs:Equations} propose the DAS-PINNs approach, which samples according to a Deep Adaptive Sampling (DAS) method and uses KRnet \cite{Tang2020DeepMapping}, a deep generative model, to approximate the PDF proportional to the residual. However, this approach is not straightforward to integrate into existing PINNs pipelines due to its dual-network framework.

In this work, we present a collocation point resampling method that scales to higher dimensions more efficiently without introducing significant computational overhead while maintaining the accuracy improvements achieved by previous approaches. We propose the Point Adaptive Collocation Method for Artificial Neural Networks (PACMANN), which uses the gradient of the squared residual to move collocation points toward areas with a large residual. In this approach, collocation point resampling is formulated as a maximization problem of the squared residual. First, the PINN is trained on a static grid of collocation points. After a certain number of iterations, this process is paused and the gradient of the squared residual is determined for the input coordinates of each collocation point. Based on the magnitude and direction provided by these gradients, points are moved to maximize the squared residual using established optimization methods. Since the residual landscape is static while training is paused, the process of moving points may be repeated iteratively. PACMANN includes four main hyperparameters: the resampling period, the optimizer for moving the collocation points, the stepsize, and the number of steps taken by the optimization algorithm. Key features of the method include its low computational cost and simplicity of integration in existing physics-informed neural network pipelines. Our approach builds on the work of Wang et al.~\cite{Wang2022IsNetwork}, who, independently of the aforementioned developments, found that iteratively updating the placement of collocation points by applying gradient ascent over the $L^{\infty}$ physics-informed loss results in a greater prediction accuracy for the Hamilton-Jacobi-Bellman equation.

First, we investigate the performance of PACMANN in combination with a variety of optimization algorithms for two low-dimensional problems: the one-dimensional Burgers' and Allen-Cahn equations. We then perform sensitivity studies on the number of collocation points and the method's hyperparameters. In addition, we demonstrate the suitability of PACMANN for high-dimensional and inverse problems. As test cases, we consider an inverse problem based on the two-dimensional Navier-Stokes equation, the Poisson's equation in five dimensions, and the three-dimensional Navier-Stokes equation. Finally, we apply PACMANN to a problem involving a re-entrant corner in a disk, showing the effectiveness of our approach on irregular domains. For all problems under consideration, we compare the performance of our method in terms of prediction accuracy and computational cost to state-of-the-art adaptive and non-adaptive sampling methods. Notably, our results show that our method matches the performance of state-of-the-art methods in terms of the accuracy/efficiency tradeoff for low-dimensional problems while efficiently scaling to high-dimensional problems, where it outperforms state-of-the-art methods.

This paper is organized as follows: In~\Cref{sec:Methodology}, we briefly review the PINNs framework, followed by a description of PACMANN. Next, in experimental studies in~\Cref{sec:Results}, we compare the accuracy and computational cost of PACMANN to other state-of-the-art sampling methods for five forward problems and an inverse problem. Finally, in ~\Cref{sec:Conclusions}, we summarize our findings.

\section{Methodology}
\label{sec:Methodology}

This section begins with a brief review of PINNs based on the framework presented by Raissi et al.~\cite{Raissi2019Physics-informedEquations} in 2019. Next, we propose the novel PACMANN.

\subsection{Physics-Informed Neural Networks (PINNs)}
PINNs approximate the solution of PDEs using neural networks. Generally, we consider PDEs of the form: find $\boldsymbol{u}$ such that
\begin{equation}
    \boldsymbol{u_t} + \mathcal{N}[\boldsymbol{u}] = 0, \quad \boldsymbol{x} \in \boldsymbol{\Omega}, \quad t \in [0,T], %\nonumber
    \label{eq:General_PDE_form}
\end{equation}
\noindent with the initial and boundary conditions
\begin{equation*}
    \begin{aligned}
        \boldsymbol{u} \left( \boldsymbol{x}, 0 \right) & = h \left( \boldsymbol{x} \right), & \quad & \boldsymbol{x} \in \boldsymbol{\Omega},\\
    \mathcal{B}[\boldsymbol{u}]\left( \boldsymbol{x},t \right) & = 0, & & \boldsymbol{x} \in \partial \boldsymbol{\Omega}, \quad t \in [0,T],
    \end{aligned}
\end{equation*}
where $\mathcal{N}[\cdot]$ is a linear or nonlinear differential operator, and $\mathcal{B}[\cdot]$ is a boundary operator corresponding to a set of boundary conditions. In addition, $\boldsymbol{x} \in \boldsymbol{\Omega} \subset \mathbb{R}^{d}$ and $t \in [0,T]$ denote the spatial and temporal coordinates, respectively, and we write $\partial \boldsymbol{\Omega}$ for the boundary of $\boldsymbol{\Omega}$. We denote the space-time domain by $\mathcal{D} \coloneqq \overline{\boldsymbol{\Omega}} \times [0,T]$.

The PINN consists of a (deep) neural network with the coordinates $(\boldsymbol{x}, t)$ as inputs and $\hat{\boldsymbol{u}}(\boldsymbol{x}, t, \boldsymbol{\theta})$ as output, approximating $\boldsymbol{u}(\boldsymbol{x}, t)$. The trainable parameters $\boldsymbol{\theta}$ of this neural network are trained by minimizing a specific loss function $\mathcal{L}(\boldsymbol{\theta})$:
\begin{equation} \label{eq:minimization_problem}
    \boldsymbol{\theta}^* = \argmin_{\boldsymbol{\theta}} \mathcal{L}\left(\boldsymbol{\theta}\right).
\end{equation}
The loss function is defined as
\begin{equation}
    \mathcal{L} \left( \boldsymbol{\theta} \right) = \lambda_{r} \mathcal{L}_{r} \left( \boldsymbol{\theta} \right) + \lambda_{ic} \mathcal{L}_{ic} \left( \boldsymbol{\theta} \right) + \lambda_{bc} \mathcal{L}_{bc} \left( \boldsymbol{\theta} \right),
    \label{eq:General_loss_function}
\end{equation}
where
\begin{align}
    & \mathcal{L}_r \left( \boldsymbol{\theta} \right) = \frac{1}{N_r} \sum_{i=1}^{N_r} \left( \hat{\boldsymbol{u_t}} \left( \boldsymbol{x}_{r}^i, t_{r}^i, \boldsymbol{\theta} \right) + \mathcal{N}[\hat{\boldsymbol{u}}] \left( \boldsymbol{x}_{r}^i, t_{r}^i, \boldsymbol{\theta} \right) \right)^2, \label{eq:Physics_loss_function} \\
    & \mathcal{L}_{ic} \left( \boldsymbol{\theta} \right) = \frac{1}{N_{ic}} \sum_{i=1}^{N_{ic}} \left( \hat{\boldsymbol{u}} \left( \boldsymbol{x}_{ic}^i, 0, \boldsymbol{\theta} \right) - h \left( \boldsymbol{x}_{ic}^i \right) \right)^2, \label{eq:Initial_loss_function} \\
    & \mathcal{L}_{bc} \left( \boldsymbol{\theta} \right) = \frac{1}{N_{bc}} \sum_{i=1}^{N_{bc}} \left( \mathcal{B}[\hat{\boldsymbol{u}}] \left( \boldsymbol{x}_{bc}^i, t_{bc}^i, \boldsymbol{\theta} \right) \right)^2 \label{eq:Boundary_loss_function}
\end{align}
represent the loss terms for the PDE residual, initial conditions, and the boundary conditions, respectively. Furthermore, $N_r$, $N_{ic}$, and $N_{bc}$ denote the numbers of collocation points of the aforementioned terms. The hyperparameters $\lambda_{r}$, $\lambda_{ic}$, and $\lambda_{bc}$ are scalar weights used to balance the loss function. Each loss term is evaluated over a set of data points, where $ \{ (\boldsymbol{x}_{r}^i, t_{r}^i) \}_{i=1}^{N_{r}} $ is a set of collocation points located in the interior of the domain, $ \{ (\boldsymbol{x}_{ic}^i,0) \}_{i=1}^{N_{ic}} $ is a set of points sampled at the initial time, and $ \{ (\boldsymbol{x}_{bc}^i, t_{bc}^i) \}_{i=1}^{N_{bc}} $ is a set sampled along the boundary. These points may be fixed during training \cite{Wu2023ANetworks}, resampled through periodic random resampling \cite{Wu2023ANetworks}, or resampled using adaptive sampling methods based on guiding information, such as the PDE residual \cite{Lu2021DeepXDE:Equations, Nabian2021EfficientSampling, Wu2023ANetworks}. Note that we assume sufficient regularity, existence of a strong-form solution of~\Cref{eq:General_PDE_form}, for the PINN loss function to be meaningful.

To train the model parameters $\boldsymbol{\theta}$, the gradient of the loss function with respect to the parameters is determined using back-propagation \cite{Rumelhart1986LearningErrors}. Next, the model parameters are updated with an optimization algorithm, often based on the gradient descent method, such as the Adam optimizer \cite{Kingma2015Adam:Optimization}. Similarly, the derivatives of $\hat{\boldsymbol{u}} \left( \boldsymbol{x}, t, \boldsymbol{\theta} \right)$ with respect to the input coordinates $\left( \boldsymbol{x}, t \right)$ as required by $\mathcal{N}[\cdot]$ and potentially $\mathcal{B}[\cdot]$ in~\Cref{eq:Physics_loss_function,eq:Boundary_loss_function} are computed using automatic differentiation; the initial conditions and the corresponding loss function~\Cref{eq:Initial_loss_function} may also depend on the temporal derivative of $\hat{\boldsymbol{u}}$, but we omit these cases for simplicity.

For problems that incorporate reference data during training, such as inverse problems, an additional loss term $\mathcal{L}_{ref}$ is added to the loss function described by~\Cref{eq:General_loss_function}, where
\begin{equation}
    \mathcal{L}_{ref} \left( \boldsymbol{\theta} \right) = \frac{1}{N_{ref}} \sum_{i=1}^{N_{ref}} \left( \hat{\boldsymbol{u}} \left( \boldsymbol{x}^i, t^i, \boldsymbol{\theta} \right) - \boldsymbol{u}_{ref} \left( \boldsymbol{x}^i, t^i \right) \right) ^2. \nonumber
\end{equation}
This term corresponds to the mean squared error between the (noisy) observed data $\boldsymbol{u}_{ref}$ at the set of data points $ \{ (\boldsymbol{x}^i, t^i) \}_{i=1}^{N_{ref}}$ and the approximation $\hat{\boldsymbol{u}} \left( \boldsymbol{x}^i, t^i, \boldsymbol{\theta} \right)$ given by the neural network. Adding this term leads to the following loss function
\begin{equation}
    \mathcal{L} \left( \boldsymbol{\theta} \right) = \lambda_{r} \mathcal{L}_{r} \left( \boldsymbol{\theta} \right) + \lambda_{ic} \mathcal{L}_{ic} \left( \boldsymbol{\theta} \right) + \lambda_{bc} \mathcal{L}_{bc} \left( \boldsymbol{\theta} \right) + \lambda_{ref} \mathcal{L}_{ref} \left( \boldsymbol{\theta} \right), \nonumber
    \label{eq:General_loss_function_ref_data_problem}
\end{equation}
where $\lambda_{ref}$ is the additional scalar weight assigned to the reference data loss term.

\subsection{Point Adaptive Collocation Method for Artificial Neural Networks (PACMANN)}

\begin{algorithm}[t]
\caption{\textbf{PACMANN with a given optimization algorithm, $P$, $s$, and $T$}}
\label{alg:Point_Marching_Adaptive_Collocation_Method}
% 1
Sample a set $\boldsymbol{X_r}$ of $N_r$ collocation points $ \{ \boldsymbol{x}_{r}^i, t_{r}^i \}_{i=1}^{N_{r}} $ using a uniform sampling method\;

% 2
\Repeat{the total number of iterations reaches the limit}{

% 3
Train the PINN for $P$ iterations\;

% 4
Determine $r^2(\boldsymbol{x},t) = \left( \boldsymbol{u_t} ( \boldsymbol{x}_{r}^i, t_{r}^i, \boldsymbol{\theta} ) + \mathcal{N}[\boldsymbol{u}] ( \boldsymbol{x}_{r}^i, t_{r}^i, \boldsymbol{\theta} ) \right)^2$, the squared PDE residual over $\boldsymbol{X_r}$\;

% 5
Find $\nabla_{\boldsymbol{x}} r^2(\boldsymbol{x},t)$ and $\frac{\partial}{\partial t} r^2(\boldsymbol{x}, t)$, the gradients of the squared PDE residual with respect to the input coordinates of the points in $\boldsymbol{X}_r$\;

% 6
Iteratively move the points in $\boldsymbol{X}_r$ according to the chosen optimization algorithm with stepsize $s$ and number of steps $T$\;

% 7
Replace points in $\boldsymbol{X}_r$ outside the domain with points sampled according to a uniform probability distribution\;

}
\end{algorithm}

\begin{figure}[t]
    \centering
    \includegraphics[width=0.45\linewidth]{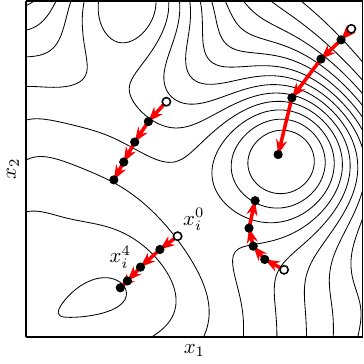}
    \caption{A schematic of PACMANN with four steps of gradient ascent on a contour plot of the squared residual.}
    \label{fig:Point_Marching_Adaptive_Collocation_Method_visual_description}
\end{figure}

In this work, we propose the \emph{Point Adaptive Collocation Method for Artificial Neural Networks}, which uses the gradient of the squared residual as guiding information to gradually move collocation points toward areas of high residuals; see~\Cref{alg:Point_Marching_Adaptive_Collocation_Method} and~\Cref{fig:Point_Marching_Adaptive_Collocation_Method_visual_description}. In particular, instead of the minimization problem in~\Cref{eq:minimization_problem}, we consider the following min-max problem for training the model parameters $\boldsymbol{\theta}$:
\begin{equation} \label{eq:min_max}
    \boldsymbol{\theta}^* = \argmin_{\boldsymbol{\theta}} \left[ \lambda_{ic} \mathcal{L}_{ic} \left( \boldsymbol{\theta} \right) + \lambda_{bc} \mathcal{L}_{bc} \left( \boldsymbol{\theta} \right) + \lambda_{ref} \mathcal{L}_{ref} \left( \boldsymbol{\theta} \right) + \lambda_{r} \max_{\boldsymbol{X_r} \subset \mathcal{D}} \mathcal{L}_{r} \left( \boldsymbol{X_r}, \boldsymbol{\theta} \right) \right].
\end{equation}
Here, only the collocation points $X_r = \{ (\boldsymbol{x}_{r}^i, t_{r}^i) \}_{i=1}^{N_{r}}$ are moved, while points such as those sampled along the boundary are fixed in place throughout training. This approach ensures that a sufficient number of points are placed along the initial and boundary conditions to accurately compute the respective loss terms.

\begin{remark}
    The min-max formulation for the residual loss $\mathcal{L}_{r}$ in~\Cref{eq:min_max} can also be interpreted as optimizing the supremum norm in $\mathcal{D}$. In our numerical experiments, the neural network function will be smooth, due to the use of the hyperbolic activation function, such that the maximum exists in $\mathcal{D}$. Then, we obviously have that
\begin{equation*}
    \begin{aligned}
        \left\|
        \hat{\boldsymbol{u_t}} + \mathcal{N}[\hat{\boldsymbol{u}}]
        \right\|_{L^\infty (\mathcal{D})}^2
        & =
        \max_{ \left( \boldsymbol{x}_r, t_r \right) \in \mathcal{D}}
        \left( \hat{\boldsymbol{u_t}} \left( \boldsymbol{x}_r, t_r \right) + \mathcal{N}[\hat{\boldsymbol{u}}] \left( \boldsymbol{x}_r, t_r \right) \right)^2
        =
        \frac{1}{N_r}
        \sum_{i=1}^{N_r}
        \max_{ \left( \boldsymbol{x}_r, t_r \right) \in \mathcal{D}}
        \left(
            \hat{\boldsymbol{u_t}} \left( \boldsymbol{x}_r, t_r \right) + \mathcal{N}[\hat{\boldsymbol{u}}] \left( \boldsymbol{x}_r, t_r \right)
        \right)^2 \\
        & =
        \max_{ \boldsymbol{X}_r \subset \mathcal{D} }
        \left[
            \frac{1}{N_r} \sum_{i=1}^{N_r}
            \left( \hat{\boldsymbol{u_t}} \left( \boldsymbol{x}_r^i, t_r^i \right) + \mathcal{N}[\hat{\boldsymbol{u}}] \left( \boldsymbol{x}_r^i, t_r^i \right) \right)^2
        \right]
        =
        \max_{\boldsymbol{X}_r \subset \mathcal{D}} \mathcal{L}_{r} \left( \boldsymbol{X_r}, \boldsymbol{\theta} \right),
    \end{aligned}
\end{equation*}
where we have omitted the dependence on the neural network parameters $\theta$ for the sake of brevity. In practice, since we never actually attain the maximum, the loss will be significantly lower than the supremum norm, due to the averaging across $\mathcal{D}$.
\end{remark}

In the PACMANN algorithm, $N_r$ collocation points are first sampled using a uniform sampling method, such as an equispaced uniform grid or the Hammersley sequence \cite[pp.~31--36]{Hammersley1964MonteMethods}. The PINN is then trained on this set of collocation points for a number of $P$ iterations. The number $P$ is a hyperparameter of the method determining the \emph{resampling period}, that is, the period after which the collocation points are resampled. Next, the training iteration is paused and the gradients of the squared residual $r^2(\boldsymbol{x},t)$ with respect to the input coordinates, given by $\nabla_{\boldsymbol{x}} r^2(\boldsymbol{x},t)$ and $\frac{\partial}{\partial t} r^2(\boldsymbol{x}, t)$, are determined for each collocation point $\left( \boldsymbol{x}, t \right) \in \boldsymbol{X}_r$. The collocation points are subsequently moved in the direction of increasing residual based on
\begin{equation}
\label{eq:Gradient_ascent_equations}
\left\{
\begin{aligned}
\boldsymbol{x}_r^{i+1} & = \boldsymbol{x}_r^i + s \nabla_{\boldsymbol{x}} r^2(\boldsymbol{x}_r^i,t_r^i), \\
t_r^{i+1} & = t_r^i + s \frac{\partial}{\partial t} r^2(\boldsymbol{x}_r^i,t_r^i). \\
\end{aligned} \right.
\end{equation}
Here, $s$ is a hyperparameter that determines the stepsize of each move; in the machine learning community, this parameter is also often called a \emph{learning rate}. Since the neural network parameters $\theta$ are kept constant during the iteration in~\Cref{eq:Gradient_ascent_equations}, the residual landscape is static and the residual gradient can be determined again for the new location of the collocation point, allowing the process to be repeated. Therefore, collocation points can be moved several times, given by the hyperparameter \emph{number of steps} $T$. If a point moves outside the domain, it is removed from the set $\boldsymbol{X}_r$ and a replacement point is sampled in the domain based on a uniform probability distribution. Next, the neural network is trained for another $P$ iterations, after which the process of moving collocation points for $T$ iterations is repeated.

In~\Cref{eq:Gradient_ascent_equations}, points are moved directly based on the gradient vector. This approach is essentially equal to applying gradient ascent. However, other gradient-based optimization algorithms can be applied. The optimization algorithms considered for PACMANN in~\Cref{sec:Results} are listed as follows. We also provide the iteration rule for an arbitrary variable $x$ and a function $f(x)$ to be maximized. Other optimization algorithms are also applicable, but we focus on the following algorithms in our numerical experiments in~\Cref{sec:Results} and~\ref{sec:Appendix_A}.
\begin{enumerate}
    \item \textbf{Gradient ascent:} The collocation points are directly moved in the direction of steepest ascent, that is, in the direction of the gradient. This algorithm updates variables using the formula:
    \begin{equation}
        x_{i+1} = x_i + s f'(x_i). \nonumber
    \end{equation}

    \vspace{10pt}

    \item \textbf{Nonlinear gradient ascent:} In this algorithm, we apply a nonlinear function to the gradient ascent algorithm to scale down large gradients, preventing points from taking large steps directly out of the domain. We refer to the algorithm as nonlinear gradient ascent. In this work, we use the hyperbolic tangent function, as follows:
    \begin{equation}
        x_{i+1} = x_i + s \tanh{ \left(f'(x_i) \right)}. \nonumber
    \end{equation}

    \vspace{10pt}

    \item \textbf{RMSprop:} Root Mean Square Propagation (RMSprop)~\cite{Hinton2012LectureMagnitude} adapts the stepsize by dividing the gradient by a weighted average of previous gradients; this serves to stabilize convergence. The algorithm consists of two steps. First, a parameter $S$ is updated. This parameter consists of a weighted average of previous gradients:
    \begin{equation}
        S_{i+1} = \beta S_i + (1 - \beta) \left(f'(x_i) \right)^2. \nonumber
    \end{equation}
    Next, the variable $x$ is updated using:
    \begin{equation}
        x_{i+1} = x_i + s \frac{f'(x_i) }{\sqrt{S_{i+1} + \epsilon}}. \nonumber
    \end{equation}
    To prevent large steps due to small values of $S$, a small value is added, represented by $\epsilon$.

    \vspace{10pt}

    \item \textbf{Momentum:} The momentum optimizer~\cite{Sutskever2013OnLearning} considers a weighted average of previous gradients at each iteration to prevent converging to local minima. First, the weighted average, $V$, is computed:
    \begin{equation}
        V_{i+1} = \beta V_{i} + (1-\beta) \left(f'(x_i) \right). \nonumber
    \end{equation}
    After updating $V$, the variable $x$ is updated:
    \begin{equation}
        x_{i+1} = x_i + s V_{i+1}. \nonumber
    \end{equation}

    \vspace{10pt}

    \item \textbf{Adam:} The Adaptive moments (Adam) \cite{Kingma2015Adam:Optimization} optimizer combines the concepts behind RMSprop and momentum. First, the value of $V$ is updated:
    \begin{equation}
        V_{i+1} = \beta_1 V_{i} + (1-\beta_1) f'(x_i). \nonumber
    \end{equation}
    Next, the parameter $S$ is updated:
    \begin{equation}
        S_{i+1} = \beta_2 S_i + (1-\beta_2) \left(f'(x_i) \right)^2. \nonumber
    \end{equation}
    Afterwards, an initialization bias correction is applied to $V_{i+1}$ and $S_{i+1}$:
    \begin{equation}
        \hat{V}_{i+1} = \frac{V_{i+1}}{1 - \beta_{1}^{i+1}}, \quad \hat{S}_{i+1} = \frac{S_{i+1}}{1 - \beta_{2}^{i+1}}. \nonumber
    \end{equation}
    Finally, the variable $x$ is updated as follows:
    \begin{equation}
        x_{i+1} = x_i + s \frac{\hat{V}_{i+1}}{\sqrt{\hat{S}_{i+1} + \epsilon}} \nonumber
    \end{equation}
    Here, a small regularization parameter $\epsilon$ is included to prevent large steps when $\hat{S}_{i+1}$ is small.

    \vspace{10pt}

    \item \textbf{Golden section search:} Golden section search~\cite[pp.~39--42]{Kochenderfer2019AlgorithmsOptimization} is a line search method that narrows down the search interval each iteration. By searching along the direction of steepest ascent, the multidimensional optimization problem is reduced to a one-dimensional problem. In this direction, the algorithm searches in an interval $[a_i, b_i]$. In the initial interval, $a_0$ is equal to the value $x_0$, for example. We determine $b_0$ using the stepsize and the gradient:
    \begin{equation}
        b_0 = a_0 + s f'(x_0). \nonumber
    \end{equation}
    Next, $f(x)$ is evaluated at two points, $x_{l,i}$ and $x_{r,i}$, determined using
    \begin{equation}
        x_{l,i} = a_i + \alpha (b_i - a_i), \quad x_{r,i} = a_i + \beta (b_i - a_i). \nonumber
    \end{equation}
    The name ``golden section search'' refers to the golden ratio, defined as $\phi = \frac{1 + \sqrt{5}}{2}$, which is incorporated in the values of $\alpha$ and $\beta$:
    \begin{equation}
        \alpha = 1 - \phi^{-1}, \quad \beta = \phi^{-1} \nonumber
    \end{equation}
    If $f(x_{l,i}) > f(x_{r,i})$, then the interval is shortened and shifted to the left:
    \begin{equation}
        a_{i+1} = a_i, \quad b_{i+1} = x_{r,i}, \quad x_{r,i+1} = x_{l,i}. \nonumber
    \end{equation}
    Otherwise, if $f(x_{l,i}) < f(x_{r,i})$, then the interval is shortened and shifted to the right:
    \begin{equation}
        a_{i+1} = x_{l,i}, \quad b_{i+1} = b_i, \quad x_{l,i+1} = x_{r,i}. \nonumber
    \end{equation}
    After updating the interval to $[a_{i+1}, b_{i+1}]$, the algorithm is repeated again. After the final iteration, the value of the variable $x$ is found by taking the middle point of the interval $[a_N, b_N]$:
    \begin{equation}
        x_N = \frac{a_N + b_N}{2}. \nonumber
    \end{equation}
    The golden ratio ensures that either $x_{l,i+1}$ ends up on $x_{r,i}$ or $x_{r,i+1}$ on $x_{l,i}$, depending on the direction of the interval shift. Since $f(x)$ has already been determined for $x_{l,i}$ and $x_{r,i}$ in the previous iteration, $f(x)$ does not have to be evaluated again for these variable values. As a result, for each iteration, the value of $f(x)$ only has to be evaluated once, which is beneficial in terms of the computational cost.
\end{enumerate}

\section{Results}
\label{sec:Results}
In this section, we evaluate the performance of PACMANN in terms of accuracy and computational cost across various PDE examples, including simple model problems, an inverse problem, a problem defined on an irregular-shaped domain, and two high-dimensional problems.
Furthermore, we vary the hyperparameters of PACMANN to showcase its behavior and compare our method with other collocation point sampling methods proposed in the related works~\cite{Lu2021DeepXDE:Equations,Wu2023ANetworks}.

To gather data on the prediction accuracy and computational cost, each method and set of hyperparameters is run ten times with varying random seeds. The prediction accuracy is compared based on the mean and standard deviation across the ten runs of the test error, measured using the $L_2$ relative error. The $L_2$ relative error is determined as follows:
\begin{equation}
    \label{eq:L2_error}
    \varepsilon_{L_2} \coloneqq \frac{\| \boldsymbol{u}_{ref} - \boldsymbol{u}_{pred} \|_2}{\| \boldsymbol{u}_{ref} \|_2},
\end{equation}
where $\boldsymbol{u}_{ref}$ is the reference solution, which is either an analytical or numerical solution depending on the problem, and $\boldsymbol{u}_{pred}$ is the predicted solution. When an analytical solution is available, we additionally consider the prediction accuracy using the $H^1$ semi-norm
\begin{equation} \label{eq:H1_error}
	\varepsilon_{H^1} \coloneqq \frac{| \boldsymbol{u}_{ref} - \boldsymbol{u}_{pred} |_1}{| \boldsymbol{u}_{ref} |_1},
	\quad \text{where} \quad
	| v |_1 \coloneqq \| \nabla v \|_2,
\end{equation}
to measure errors in the first derivatives of the solution. To compare $\boldsymbol{u}_{ref}$ and $\boldsymbol{u}_{pred}$ we employ an equispaced uniform grid of 10\,000 collocation points.
The mean runtime of training over the ten runs serves as an indication of the computational cost of a particular sampling method. Our code is based on the PINNs library DeepXDE~\cite{Lu2021DeepXDE:Equations} using PyTorch \cite{Paszke2019PyTorch:Library} version 1.12.1 as the backend. It is publicly available on GitHub at \url{https://github.com/CoenVisser/PACMANN}. The models were trained using NVIDIA Tesla V100S GPUs on TU Delft's high-performance computer DelftBlue~\cite{DelftHighPerformanceComputingCentreDHPC2024DelftBlue2}.

For all experiments, the training is split into five phases of 10\,000 iterations, consisting of 7000 iterations of Adam with a learning rate of $10^{-3}$ followed by 3000 iterations of the Limited-memory Broyden-Fletcher-Goldfarb-Shanno (L-BFGS) algorithm~\cite{Liu1989OnOptimization}, totaling 50\,000 iterations. Throughout training, a resampling period $P$ of 50 iterations is maintained. Furthermore, we only resample the collocation points when training the neural network parameters with Adam. Resampling while training with the L-BFGS optimizer would disrupt the convergence of the algorithm due to the change in loss landscape by evaluating the PDE loss term at different collocation points. In this work, in all experiments, the hyperbolic tangent is used as the activation function. In addition, the hyperparameter settings used for RMSprop, momentum, and Adam with PACMANN are given in~\Cref{tab:hyperparam_optimizers}. Note that we have not conducted an extensive study varying the neural network architecture. Instead, we have used the architectures provided in the corresponding test cases of DeepXDE, assuming that these were already optimized appropriately.

\begin{table}[t]
\centering
\begin{tabular}{lllll}
\toprule
\textbf{Optimizer} & $\boldsymbol{\beta}$ & $\boldsymbol{V_0}$ & $\boldsymbol{S_0}$ & $\boldsymbol{\epsilon}$ \\ \midrule
RMSprop  & 0.999                                 & - & 0 & $10^{-8}$ \\
Momentum & 0.9                                   & 0 & - & $10^{-8}$ \\
Adam     & $\beta_1 = 0.9$, $\beta_2 = 0.999$ & 0 & 0 & $10^{-8}$ \\
\bottomrule
\end{tabular}
\caption{Hyperparameter settings for the RMSprop, momentum, and Adam optimizers.}
\label{tab:hyperparam_optimizers}
\end{table}

In the numerical experiments, PACMANN with Adam consistently achieves the lowest error compared to the other optimization algorithms considered. Therefore, to preserve clarity, figures and tables which compare the various sampling methods only contain the results for Adam with our method. The figures comparing the accuracy behavior of the other optimization algorithms discussed in~\Cref{sec:Methodology} for varying numbers of collocation points and resampling periods are provided in \ref{sec:Appendix_A}.

Infrequently, the random neural network weight initialization prevents the PINN from learning the solution, which results in a test error several orders of magnitude larger than the test error obtained with other weight initializations. This has been observed for all sampling methods considered in this study and is characterized by volatile loss behavior or large static loss terms. When this occurs, the corresponding training run is repeated with a different random seed.

\subsection{1D Burgers' equation} \label{sec:Burgers_equation}

\begin{table}[t]
\centering
\begin{tabular}{@{}lccc@{}}
\toprule
\multicolumn{1}{c}{\multirow{2}{*}{Sampling method}} & \multicolumn{2}{c}{$L_2$ relative error} & \multirow{2}{*}{\begin{tabular}[c]{@{}c@{}}Mean \\ runtime {[}s{]}\end{tabular}} \\ \cmidrule(lr){2-3}
\multicolumn{1}{c}{} & Mean & 1 SD &  \\ \midrule
\multicolumn{1}{l}{Uniform grid} & 25.9\% & 14.2\% & 425 \\
\multicolumn{1}{l}{Hammersley grid} & 0.61\% & 0.53\% & 443 \\
\multicolumn{1}{l}{Random resampling} & 0.40\% & 0.35\% & \textbf{423} \\ \cmidrule(r){1-4}
\multicolumn{1}{l}{RAR} & 0.11\% & \textbf{0.05\%} & 450 \\
\multicolumn{1}{l}{RAD} & 0.16\% & 0.10\% & 463 \\
\multicolumn{1}{l}{RAR-D} & 0.24\% & 0.21\% & 503 \\ \cmidrule(r){1-4}
\multicolumn{1}{l}{PACMANN-Adam} & \textbf{0.07\%} & \textbf{0.05\%} & 461 \\ \bottomrule
\end{tabular}
\caption{Overview of the mean and standard deviation of the test error and the mean runtime for each sampling method for the Burgers' equation. The best result in each column is marked in boldface.}
\label{tab:Burgers_overview_L2_error_samplers}
\end{table}

\begin{table}[t]
\centering
\begin{tabular}{@{}lccccc@{}}
\toprule
\multicolumn{1}{c}{\multirow{2}{*}{PACMANN optimizer}} & \multicolumn{2}{c}{$L_2$ relative error} & \multirow{2}{*}{\begin{tabular}[c]{@{}c@{}}Mean \\ runtime {[}s{]}\end{tabular}} & \multicolumn{2}{c}{Hyperparameters} \\ \cmidrule(lr){2-3} \cmidrule(l){5-6}
\multicolumn{1}{c}{} & Mean & 1 SD &  & Stepsize $s$ & No. of steps $T$ \\ \midrule
\multicolumn{1}{l}{Gradient ascent} & 0.30\% & 0.17\% & \multicolumn{1}{c}{\textbf{436}} & $10^{-6}$ & 1 \\
\multicolumn{1}{l}{Nonlinear gradient ascent} & 0.10\% & 0.06\% & \multicolumn{1}{c}{453} & $10^{-4}$ & 5 \\
\multicolumn{1}{l}{RMSprop} & 0.10\% & \textbf{0.03\%} & \multicolumn{1}{c}{442} & $10^{-6}$ & 10 \\
\multicolumn{1}{l}{Momentum} & 0.18\% & 0.24\% & \multicolumn{1}{c}{448} & $10^{-6}$ & 5 \\
\multicolumn{1}{l}{Adam} & \textbf{0.07\%} & 0.05\% & \multicolumn{1}{c}{461} & $10^{-5}$ & 15 \\
\multicolumn{1}{l}{Golden section search} & 0.34\% & 0.17\% & \multicolumn{1}{c}{460} & $10^{-7}$ & 5 \\ \bottomrule
\end{tabular}
\caption{Overview of the mean and standard deviation of the test error and the mean runtime achieved by PACMANN for the Burgers' equation with the optimization methods listed in~\Cref{sec:Methodology}. The best result in each column is marked in boldface. For each optimization method, we report the stepsize and number of steps that achieve the lowest test error.}
\label{tab:Burgers_overview_L2_error_optimizers}
\end{table}

\begin{figure}[t]
     \centering
     \begin{subfigure}[b]{0.49\textwidth}
         \centering
         \includegraphics[width=\textwidth]{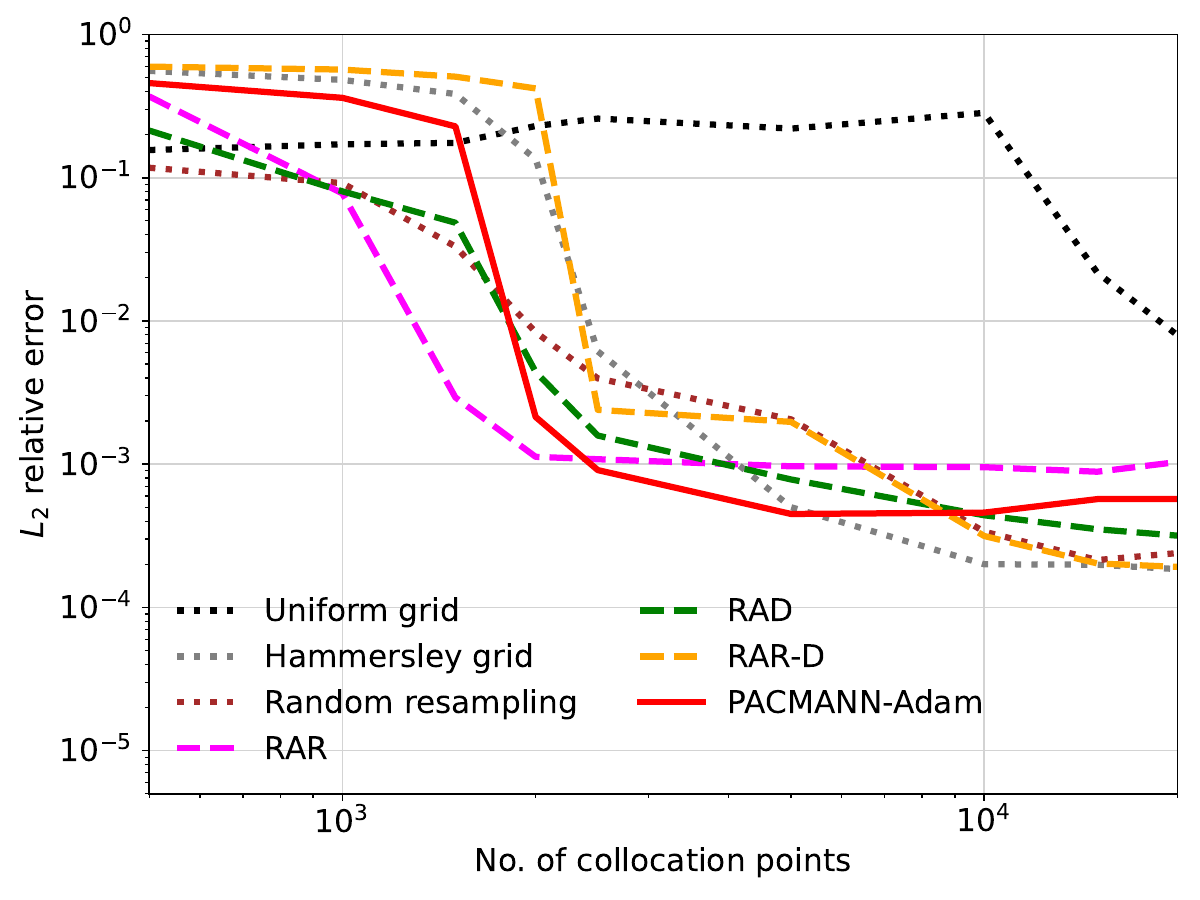}
         \caption{}
         \label{fig:Burgers_L2_vs_collocation}
     \end{subfigure}
     \begin{subfigure}[b]{0.49\textwidth}
         \centering
         \includegraphics[width=\textwidth]{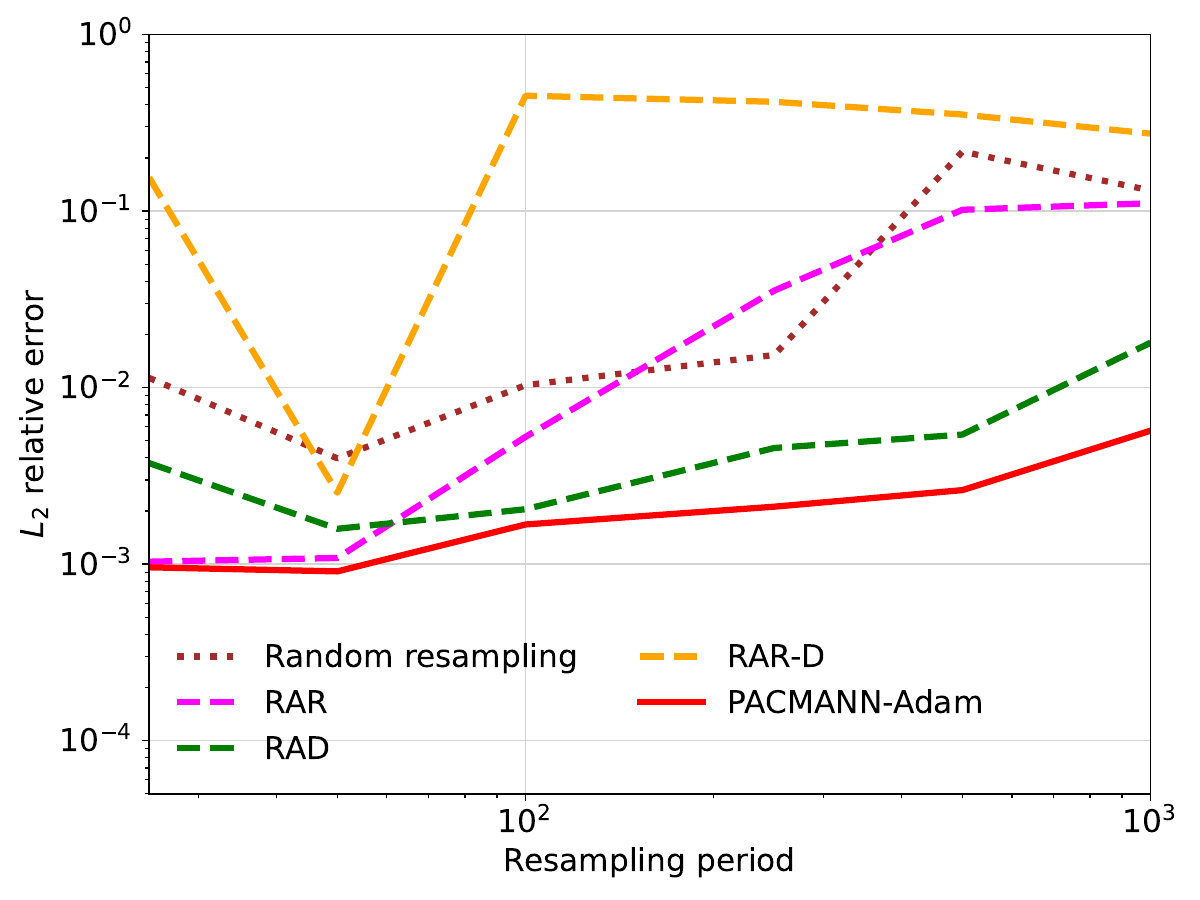}
         \caption{}
         \label{fig:Burgers_L2_vs_period}
     \end{subfigure}
        \caption{Mean of the test error for each of the sampling methods for a varying (a) number of collocation points and (b) resampling period for the Burgers' equation example.}
        \label{fig:Burgers_L2_vs_main_hyperparam}
\end{figure}

\begin{figure}[t]
     \centering
     \begin{subfigure}[b]{0.49\textwidth}
         \centering
         \includegraphics[width=\textwidth]{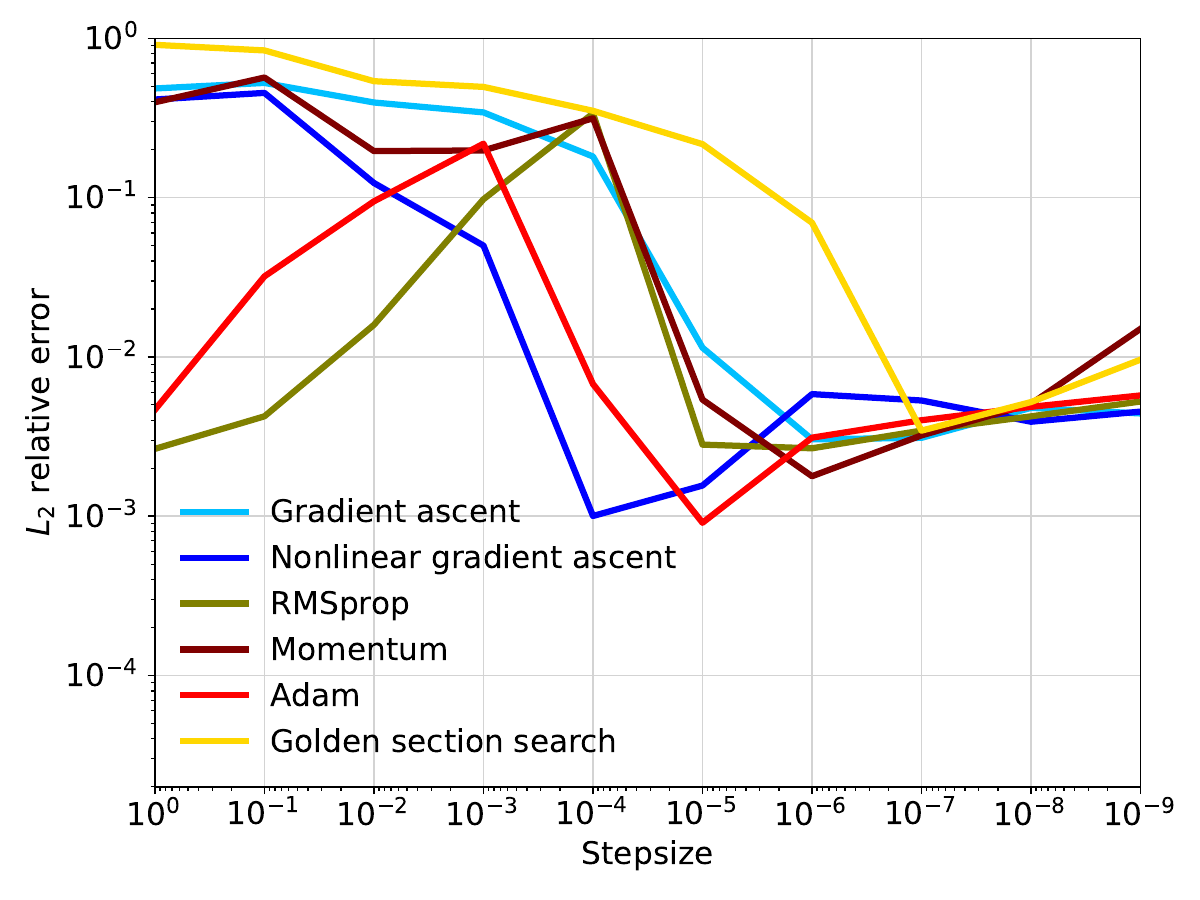}
         \caption{}
         \label{fig:Burgers_L2_vs_stepsize}
     \end{subfigure}
     \begin{subfigure}[b]{0.49\textwidth}
         \centering
         \includegraphics[width=\textwidth]{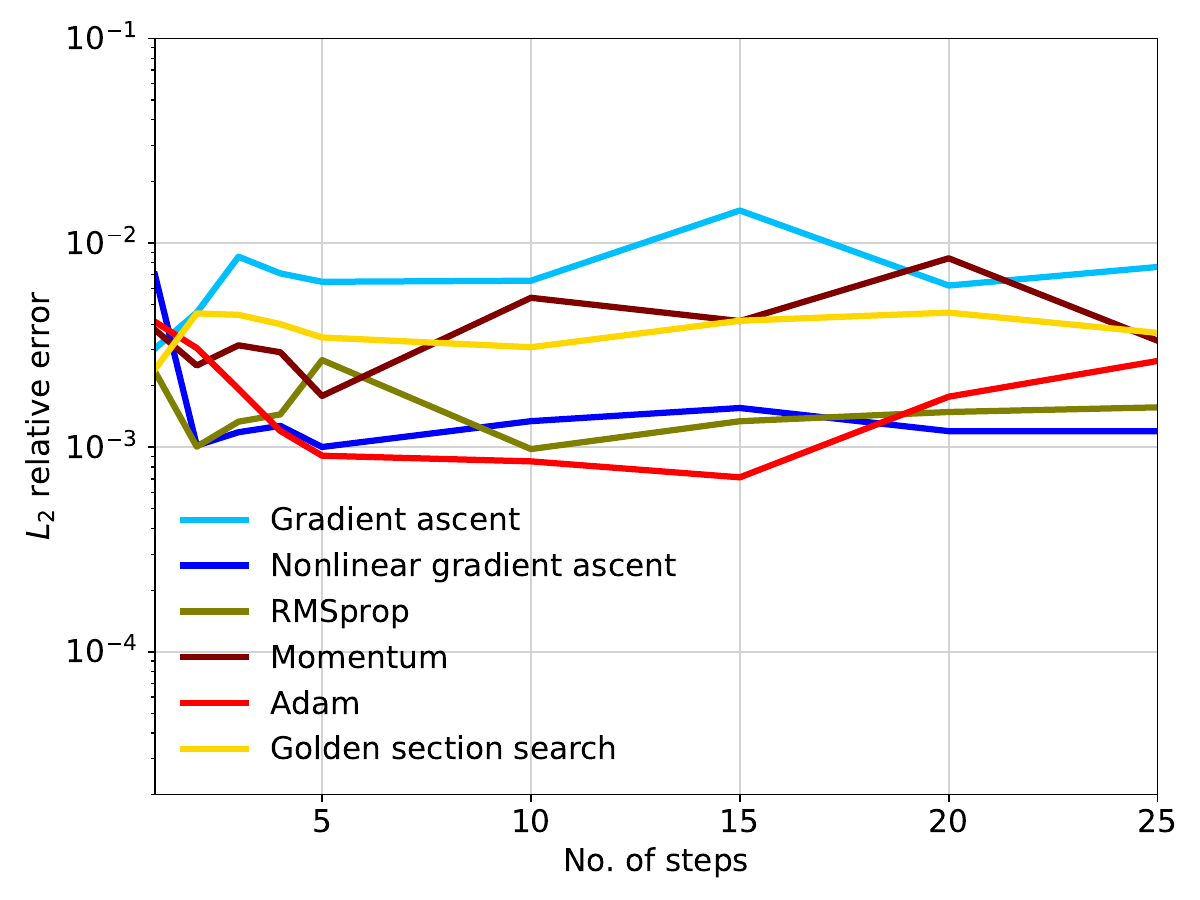}
         \caption{}
         \label{fig:Burgers_L2_vs_iter}
     \end{subfigure}
        \caption{Mean of the test error for PACMANN with the different optimization algorithms listed in~\Cref{sec:Methodology} for a varying (a) stepsize and (b) number of steps for the Burgers' equation example.}
        \label{fig:Burgers_L2_vs_step_iter}
\end{figure}

\begin{figure}[t]
     \centering
     \begin{subfigure}[b]{0.49\textwidth}
         \centering
         \includegraphics[width=\textwidth]{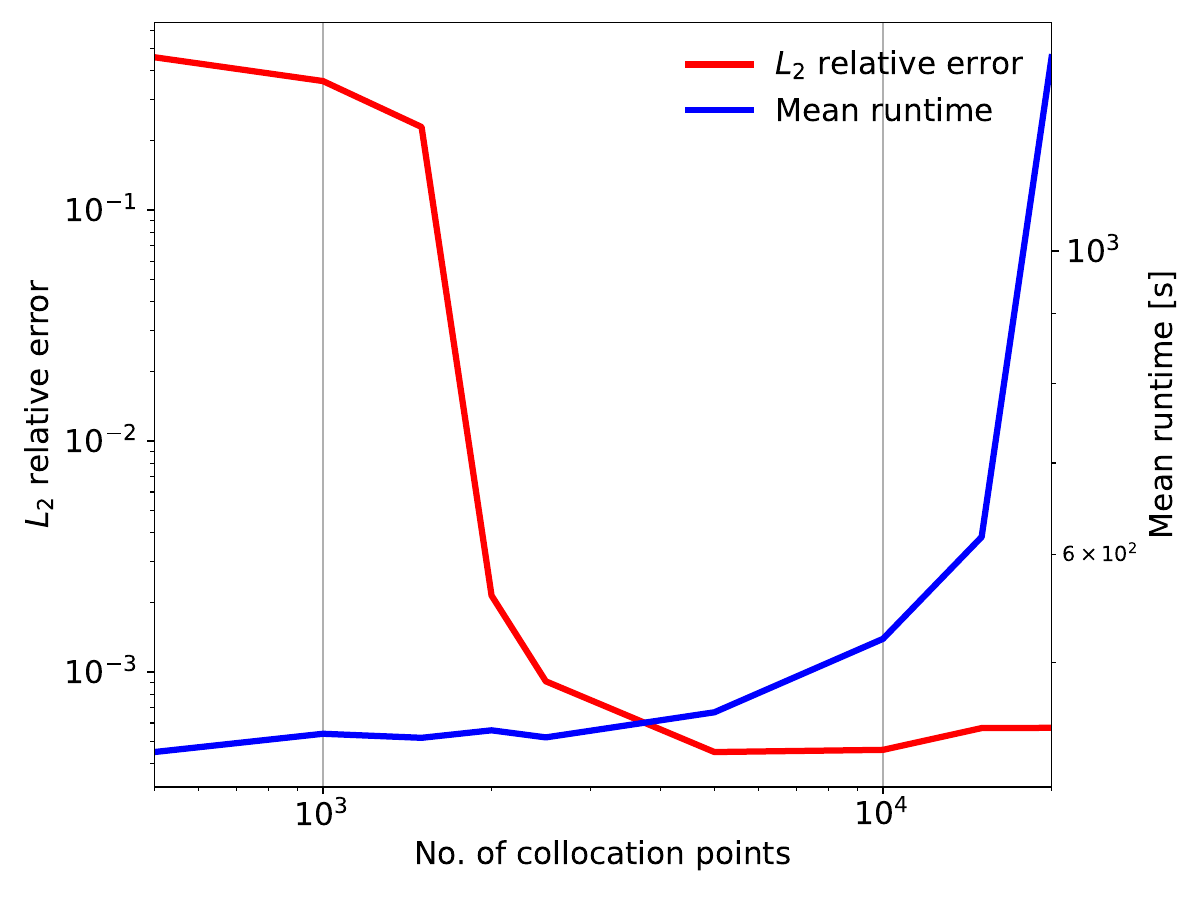}
         \caption{}
         \label{fig:Burgers_L2_Duration_vs_col}
     \end{subfigure}
     \begin{subfigure}[b]{0.49\textwidth}
         \centering
         \includegraphics[width=\textwidth]{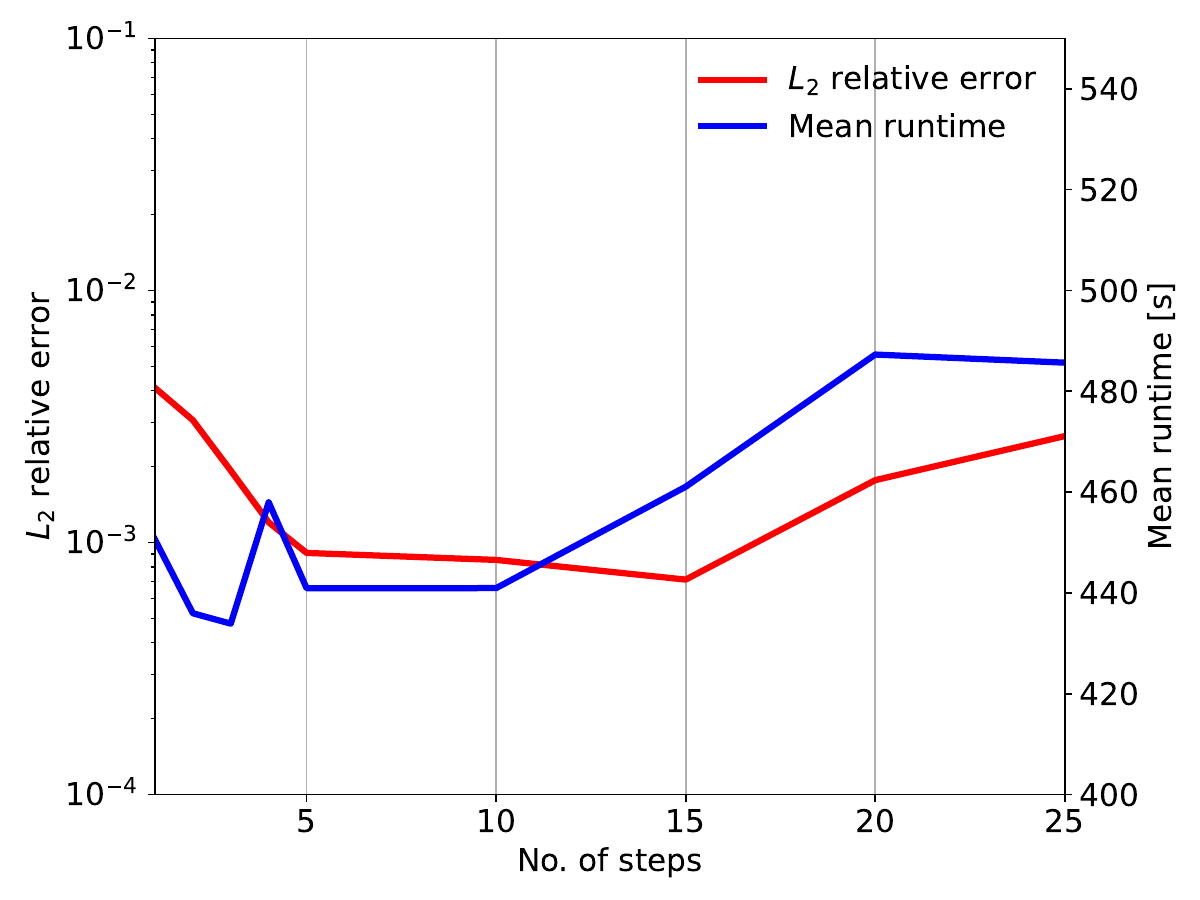}
         \caption{}
         \label{fig:Burgers_L2_Duration_vs_iter}
     \end{subfigure}
        \caption{Mean of the test error and the runtime for PACMANN with Adam for a varying (a) number of collocation points and (b) number of steps for the Burgers' equation example.}
        \label{fig:Burgers_L2_vs_Adam_col_iter}
\end{figure}

\begin{figure}[p]
     \centering
     \begin{subfigure}[b]{0.47\textwidth}
         \centering
         \includegraphics[width=\textwidth]{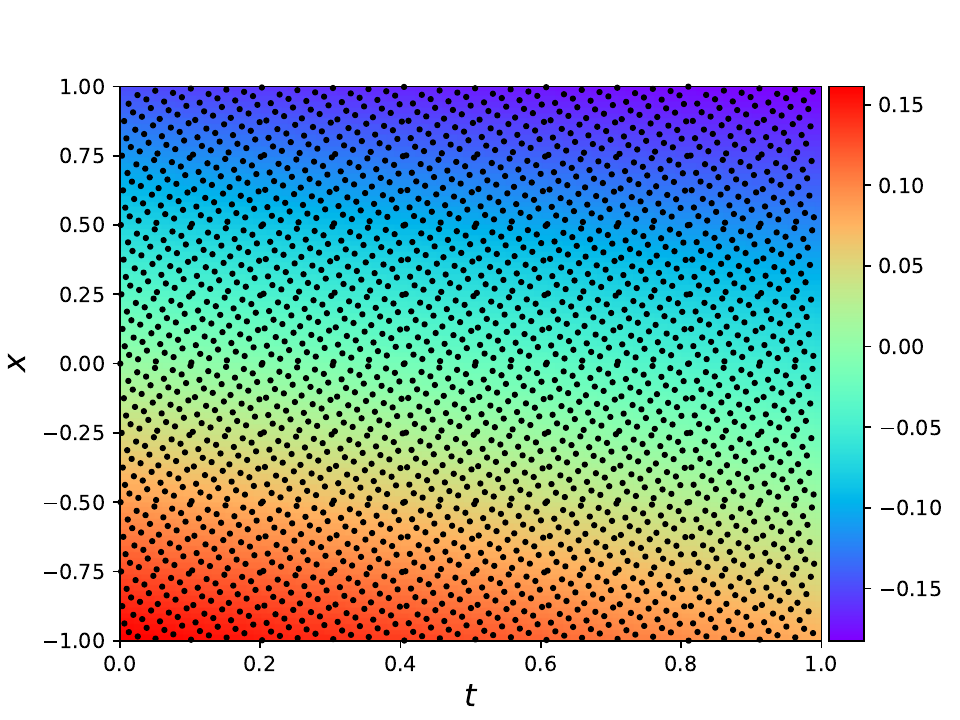}
         \caption{}
         \label{fig:Burgers_0_iter}
     \end{subfigure}
     \begin{subfigure}[b]{0.47\textwidth}
         \centering
         \includegraphics[width=\textwidth]{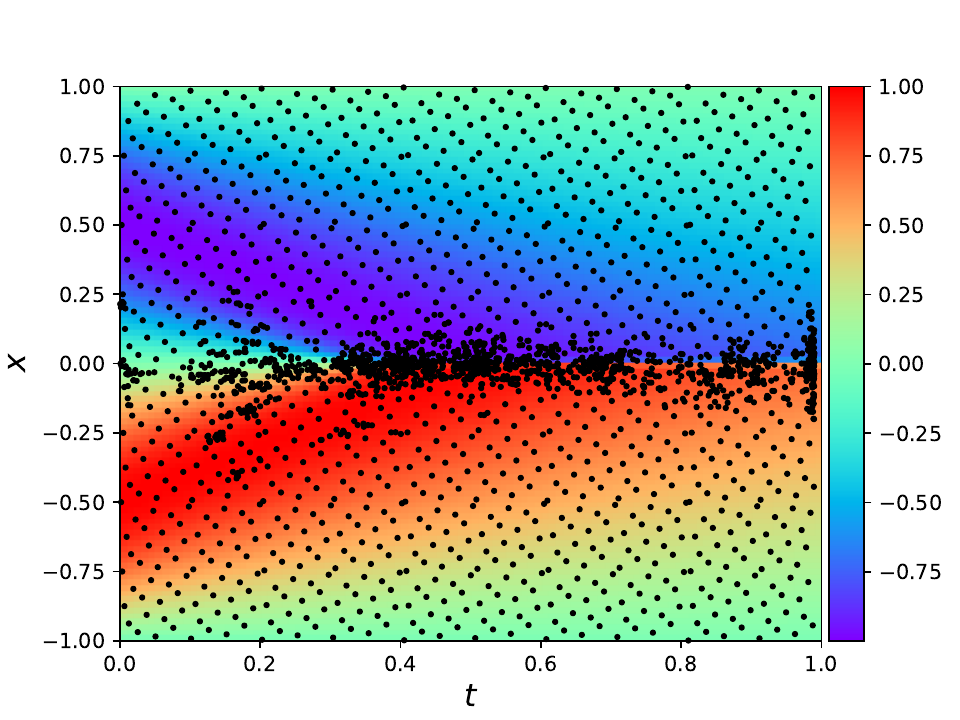}
         \caption{}
         \label{fig:RAR_burgers}
     \end{subfigure}
     \begin{subfigure}[b]{0.47\textwidth}
         \centering
         \includegraphics[width=\textwidth]{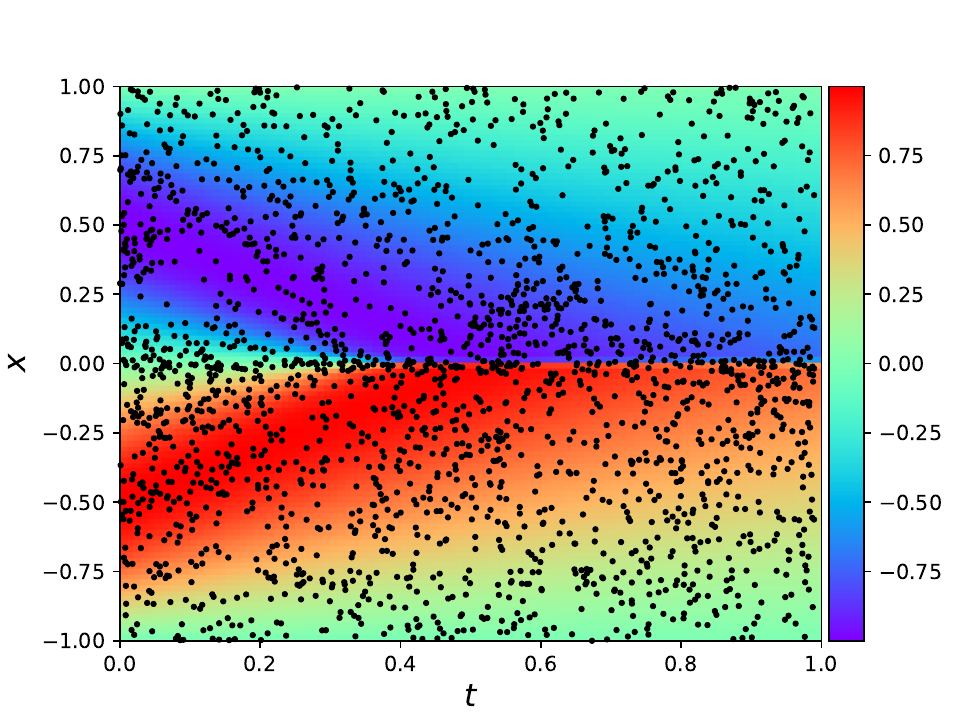}
         \caption{}
         \label{fig:RAD_burgers}
     \end{subfigure}
     \begin{subfigure}[b]{0.47\textwidth}
         \centering
         \includegraphics[width=\textwidth]{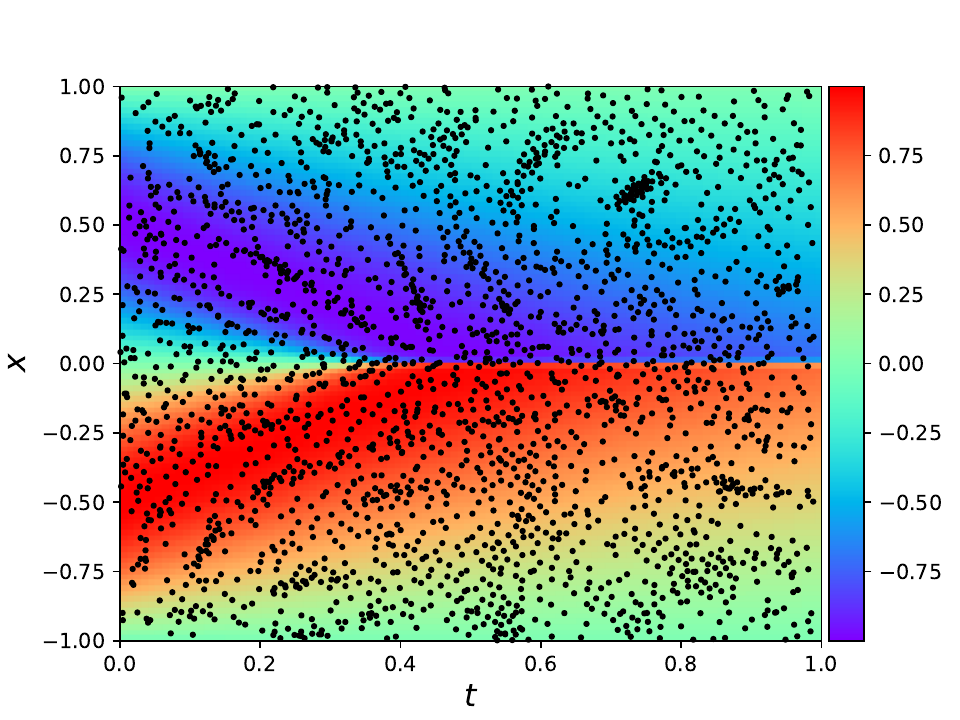}
         \caption{}
         \label{fig:Adam_burgers}
     \end{subfigure}
        \caption{Location of the 2500 collocation points (a) before training, and after training with (b) RAR, (c) RAD, and (d) PACMANN with Adam for the Burgers' equation example. The color indicates the values of the predicted solution.}
        \label{fig:Burgers_collocation_locations}
\end{figure}

\begin{figure}[p]
     \centering
     \begin{subfigure}[b]{0.47\textwidth}
         \centering
         \includegraphics[width=\textwidth]{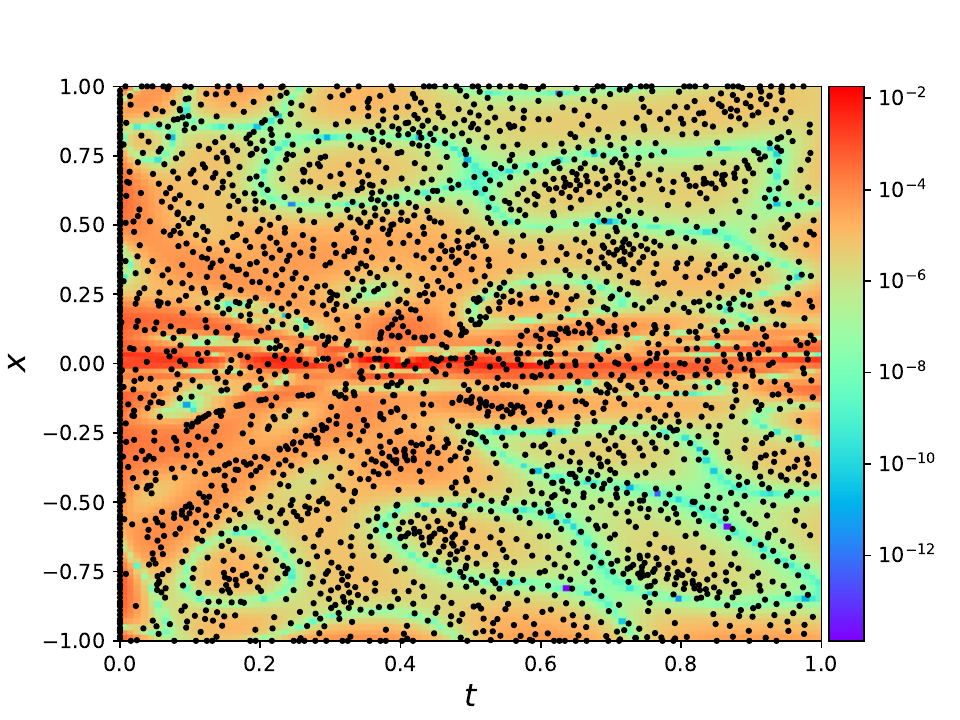}
         \caption{}
         \label{fig:R2_Mid_Training}
     \end{subfigure}
     \begin{subfigure}[b]{0.47\textwidth}
         \centering
         \includegraphics[width=\textwidth]{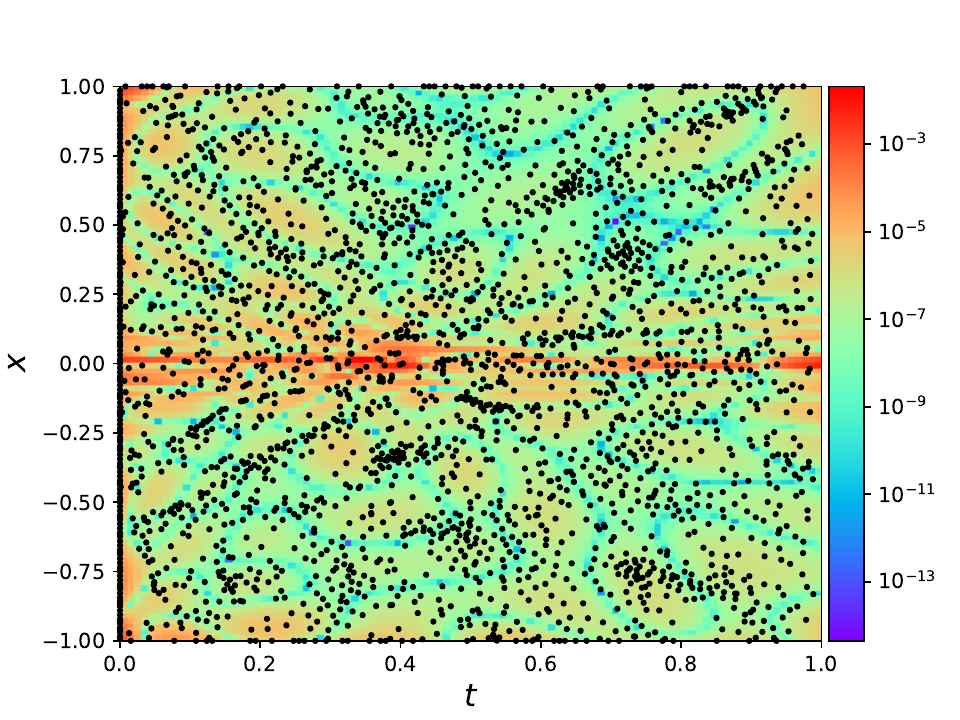}
         \caption{}
         \label{fig:R2_Post_Training}
     \end{subfigure}
        \caption{Location of the collocation points (a) mid-training at 25\,000 iterations, and (b) after training using PACMANN with Adam for the Burgers' equation example. The color indicates the values of the squared residual.}
        \label{fig:Burgers_squared_residual_plots}
\end{figure}

We first consider the one-dimensional Burgers' equation:
\begin{equation}
\label{eq:Burgers_equation_problem_statement}
\left\{
\begin{aligned}
 u_t + u u_x & = \nu u_{xx}, & \quad & x \in [-1, 1], \quad t \in [0,1], \\
 u(x,0) & = -\sin(\pi x), & \\
 u(-1,t) & = u(1,t) = 0. & \\
\end{aligned} \right.
\end{equation}
Here, $\nu$ is the diffusion coefficient or kinematic viscosity, set at $\nu = 0.01/\pi$. For this problem, we employ 2500 collocation points, 80 boundary points, and 160 initial condition points. The neural network architecture used for this example consists of four hidden layers of 64 neurons. To compute the $L_2$ relative error, we compare the network prediction to a numerical solution generated using a spectral solver by Raissi et al.~\cite{Raissi2019Physics-informedEquations}.

The mean and standard deviation of the test error and the mean runtime for each sampling method are given in~\Cref{tab:Burgers_overview_L2_error_samplers}. Out of the non-adaptive sampling methods tested, only the static uniform grid is unable to capture the solution satisfactorily, resulting in a high test error of $25.9\%$. The other non-adaptive methods, the static Hammersley grid and random resampling, attain a significantly lower test error. Overall, our method, in combination with the Adam optimizer and a stepsize of $10^{-5}$, achieves the lowest test error. It achieves a lower error than the next-best sampling method, RAR, at a slightly higher computational cost.

\Cref{tab:Burgers_overview_L2_error_optimizers} compares the performance of the various optimization algorithms for PACMANN in terms of accuracy and efficiency. We note that the nonlinear gradient ascent and the RMSprop optimizers with our proposed method both achieve a competitive test error and computational cost compared to the second best approach in~\Cref{tab:Burgers_overview_L2_error_samplers}, that is, RAR.

Next, we test the behavior of the different sampling methods by varying the number of collocation points from 500 to 20\,000; see~\Cref{fig:Burgers_L2_vs_collocation}. We observe that RAR initially reduces the test error the fastest, but plateaus at a higher error than the other methods under consideration. PACMANN in combination with the Adam optimizer significantly improves the prediction accuracy as the number of collocation points is increased from 1500 to 2000, after which it slowly increases the prediction accuracy. Only at a large number of collocation points (10\,000) are the adaptive sampling methods slightly outperformed by the static Hammersley grid.

\Cref{fig:Burgers_L2_vs_period} depicts the behavior of the various sampling methods as the resampling period is increased from 25 to 1000. Generally, most sampling methods lose accuracy as the period is increased. Notably, RAR-D performs significantly better at a period of 50 iterations compared to other resampling periods. PACMANN performs best for all resampling periods considered, losing accuracy slower than the other sampling methods.

In the following, we compare the behavior of the optimizers listed in~\Cref{sec:Methodology} when varying the hyperparameters stepsize $s$ and number of steps $T$. First, in~\Cref{fig:Burgers_L2_vs_stepsize}, we test the accuracy of the optimizers for different stepsizes $s$ ranging from $1$ to $10^{-9}$. This figure demonstrates that the stepsize has a significant influence on the prediction accuracy achieved by PACMANN. We note that, depending on the optimizer used, a different stepsize is optimal, such as $10^{-6}$ for momentum or $10^{-5}$ for Adam. Furthermore, the behavior of the optimizers at stepsizes near $1$ is split into two groups, with Adam and RMSprop gaining accuracy as the stepsize is increased, while others continue to lose accuracy. This phenomenon is explained by the number of collocation points that leave the domain while points are moved by PACMANN. At these large stepsizes, all 2500 collocation points exit the domain when using Adam or RMSprop. Since PACMANN uses a uniform probability distribution to determine the location of the replacement collocation points, the test error approaches the accuracy of random resampling ($0.40\% \pm 0.35\%$); cf.~\Cref{tab:Burgers_overview_L2_error_samplers}. In contrast, when applying the other optimization algorithms with these stepsizes, only a small portion of the collocation points exit the domain, of order $\mathcal{O}(10)$.
These few points are not sufficient for the random resampling to significantly affect the test error achieved. The difference in the number of points that exit the domain may be explained by the relatively large optimal stepsize and number of steps for RMSprop and Adam compared to the other optimization algorithms.

We test the influence of the number of steps $T$ by ranging it from 1 to 25. \Cref{fig:Burgers_L2_vs_iter} shows that certain optimizers benefit from more steps, such as nonlinear gradient ascent and Adam. Others remain at a near-constant accuracy or lose accuracy with additional steps. Based on~\Cref{fig:Burgers_L2_vs_iter}, we note that the number of steps generally has a smaller impact on the test error achieved compared to the stepsize.

Furthermore, \Cref{fig:Burgers_L2_Duration_vs_col,fig:Burgers_L2_Duration_vs_iter} depict the prediction accuracy and computational cost of PACMANN with the Adam optimizer for varying numbers of collocation points and steps. \Cref{fig:Burgers_L2_Duration_vs_col} demonstrates that increasing the number of collocation points reduces the $L_2$ relative error before reaching a plateau, beyond which the computational cost rises steeply without a further increase in accuracy. Importantly, we point out that the accuracy of Adam with five steps is nearly the same as its accuracy at 15 steps, see ~\Cref{fig:Burgers_L2_Duration_vs_iter}. Thus, we recommend taking fewer steps to save on computational cost in more complex problems.

Next, we compare visually the distribution of collocation points after training. \Cref{fig:Burgers_0_iter} shows the locations before training when the collocation points are laid out based on the Hammersley sequence. \Cref{fig:RAR_burgers,fig:RAD_burgers,fig:Adam_burgers} show the locations of the collocation points after training for RAR, RAD, and PACMANN with Adam, respectively. While RAR clusters the points at the steepest region of the solution, RAD and our method tend to create several smaller clusters.
In contrast to RAD, our method also forms clusters of points in regions with typically lower residuals, see ~\Cref{fig:R2_Mid_Training}, indicative of local maxima of the squared residual. \Cref{fig:R2_Post_Training} shows that, after training the PINN, these previously observed local maxima have reduced. For our approach, we note the similarity between the shape of the clusters and the solution itself.

\subsection{1D Allen-Cahn equation} \label{sec:Allen_Cahn_equation}

\begin{table}[t]
\centering
\begin{tabular}{@{}lccc@{}}
\toprule
\multicolumn{1}{c}{\multirow{2}{*}{Sampling method}} & \multicolumn{2}{c}{$L_2$ relative error} & \multirow{2}{*}{\begin{tabular}[c]{@{}c@{}}Mean \\ runtime {[}s{]}\end{tabular}} \\ \cmidrule(lr){2-3}
\multicolumn{1}{c}{} & Mean & 1 SD &  \\ \midrule
\multicolumn{1}{l}{Uniform grid} & 44.34\% & 18.58\% & 634 \\
\multicolumn{1}{l}{Hammersley grid} & 0.47\% & 0.26\% & \textbf{591} \\
\multicolumn{1}{l}{Random resampling} & 0.42\% & 0.28\% & 592 \\ \cmidrule(r){1-4}
\multicolumn{1}{l}{RAR} & 0.44\% & 0.27\% & 576 \\
\multicolumn{1}{l}{RAD} & 0.93\% & 0.69\% & 655 \\
\multicolumn{1}{l}{RAR-D} & 0.28\% & 0.13\% & 632 \\ \cmidrule(r){1-4}
\multicolumn{1}{l}{PACMANN-Adam} & \textbf{0.16\%} & \textbf{0.07\%} & 632 \\ \bottomrule
\end{tabular}
\caption{Overview of the mean and standard deviation of the test error and the mean runtime for each sampling method for the Allen-Cahn equation example. The best result in each column is marked in boldface.}
\label{tab:Allen_Cahn_overview_L2_error_resamplers}
\end{table}

\begin{table}[t]
\centering
\begin{tabular}{@{}lccccc@{}}
\toprule
\multicolumn{1}{c}{\multirow{2}{*}{PACMANN optimizer}} & \multicolumn{2}{c}{$L_2$ relative error} & \multirow{2}{*}{\begin{tabular}[c]{@{}c@{}}Mean \\ runtime {[}s{]}\end{tabular}} & \multicolumn{2}{c}{Hyperparameters} \\ \cmidrule(lr){2-3} \cmidrule(l){5-6}
\multicolumn{1}{c}{} & Mean & 1 SD &  & Stepsize $s$ & No. of steps $T$ \\ \midrule
\multicolumn{1}{l}{Gradient ascent} & 0.46\% & 0.24\% & \multicolumn{1}{c}{574} & $10^{-8}$ & 5 \\
\multicolumn{1}{l}{Nonlinear gradient ascent} & 0.42\% & 0.24\% & \multicolumn{1}{c}{602} & $10^{-7}$ & 5 \\
\multicolumn{1}{l}{RMSprop} & 0.29\% & 0.20\% & \multicolumn{1}{c}{595} & $10^{-6}$ & 5 \\
\multicolumn{1}{l}{Momentum} & 0.36\% & 0.17\% & \multicolumn{1}{c}{\textbf{567}} & $10^{-7}$ & 5 \\
\multicolumn{1}{l}{Adam} & \textbf{0.16\%} & \textbf{0.07\%} & \multicolumn{1}{c}{632} & $10^{-5}$ & 5 \\
\multicolumn{1}{l}{Golden section search} & 0.37\% & 0.29\% & \multicolumn{1}{c}{635} & $10^{-7}$ & 15 \\ \bottomrule
\end{tabular}
\caption{Overview of the mean and standard deviation of the test error and the mean runtime achieved by PACMANN for each optimization algorithm listed in~\Cref{sec:Methodology}. The best result in each column is marked in boldface. We also include the optimal values for the stepsize and the number of steps, for the Allen-Cahn equation example.}
\label{tab:Allen_Cahn_overview_L2_error_optimizers}
\end{table}

\begin{figure}[t]
     \centering
     \begin{subfigure}[b]{0.49\textwidth}
         \centering
         \includegraphics[width=\textwidth]{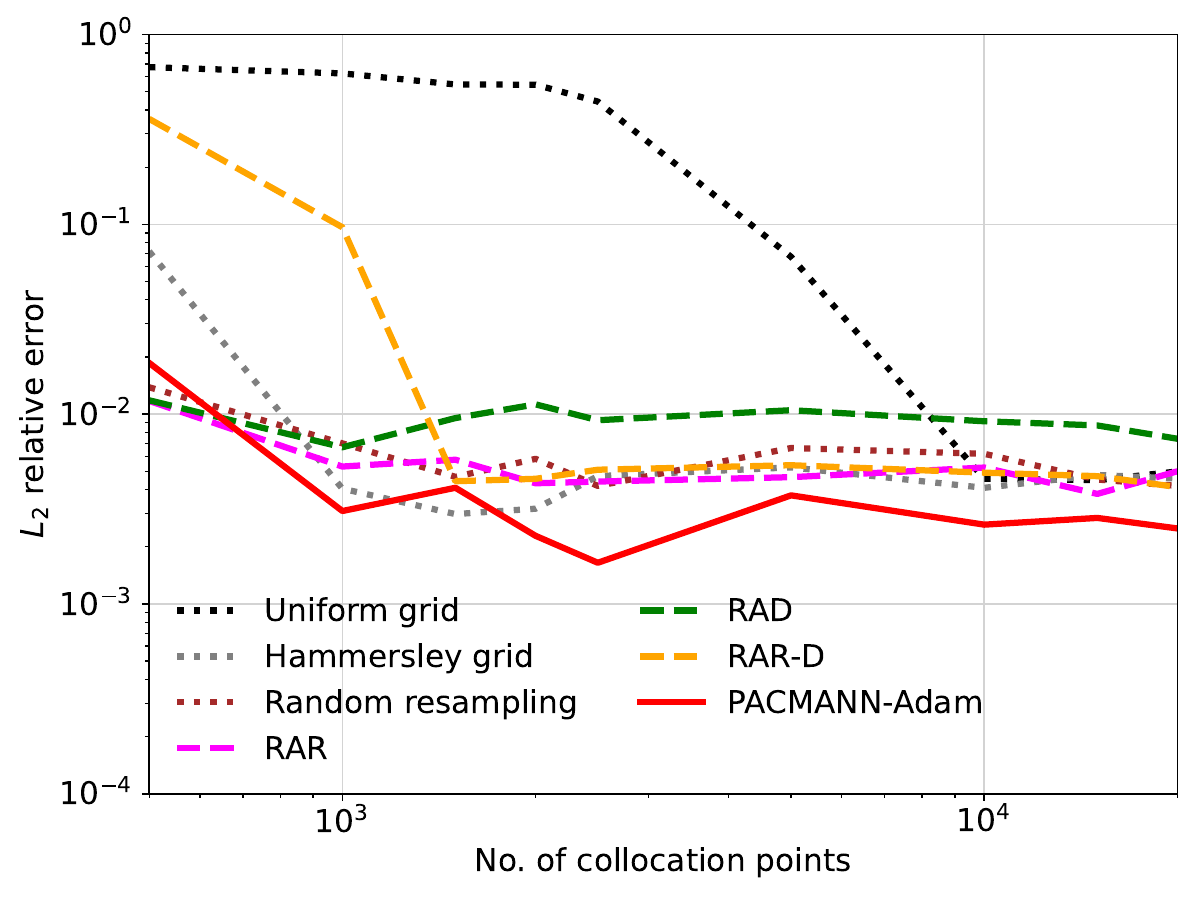}
         \caption{}
         \label{fig:Allen_Cahn_L2_vs_collocation}
     \end{subfigure}
     \begin{subfigure}[b]{0.49\textwidth}
         \centering
         \includegraphics[width=\textwidth]{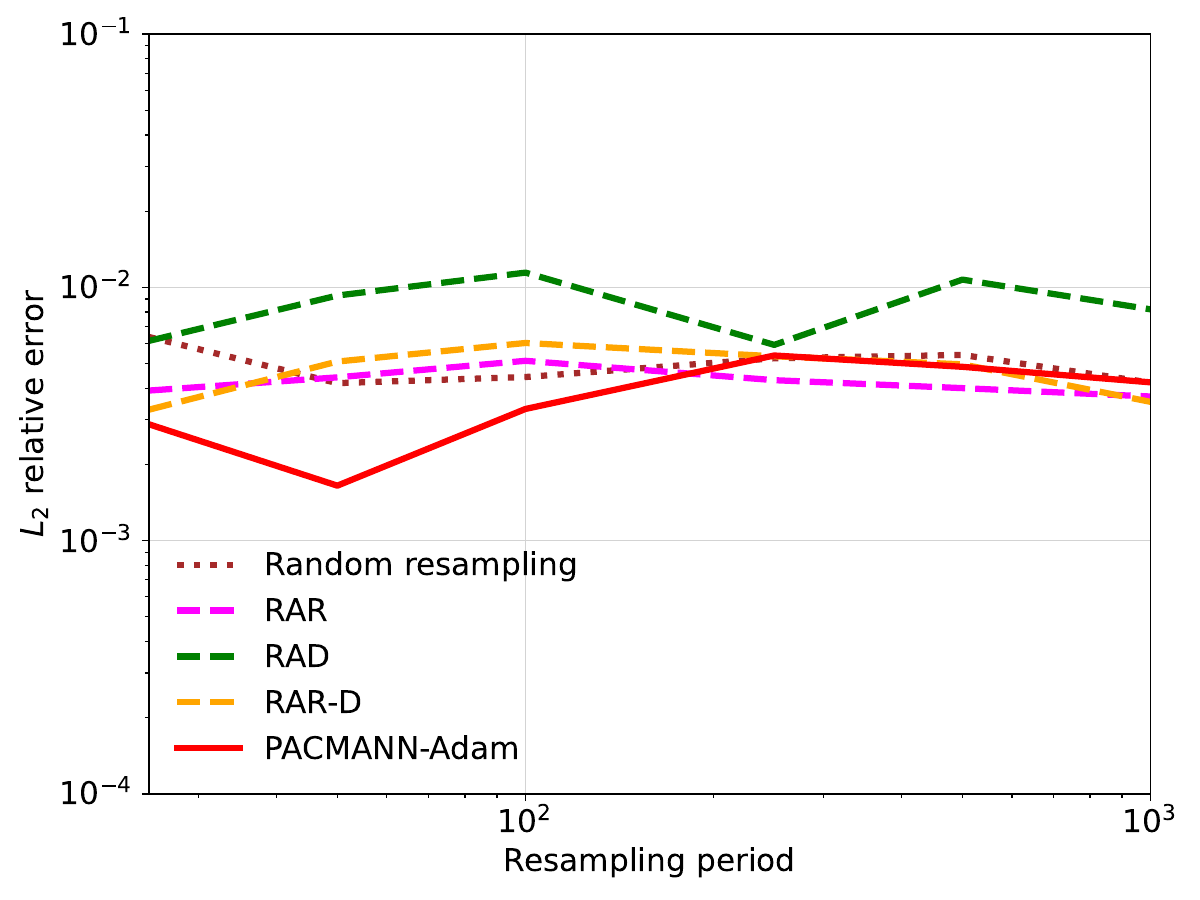}
         \caption{}
         \label{fig:Allen_Cahn_L2_vs_period}
     \end{subfigure}
        \caption{Mean of the test error for each of the sampling methods for a varying (a) number of collocation points and (b) resampling period for the Allen-Cahn equation example.}
        \label{fig:Allen_cahn_L2_vs_main_hyperparam}
\end{figure}

\begin{figure}[t]
     \centering
     \begin{subfigure}[b]{0.49\textwidth}
         \centering
         \includegraphics[width=\textwidth]{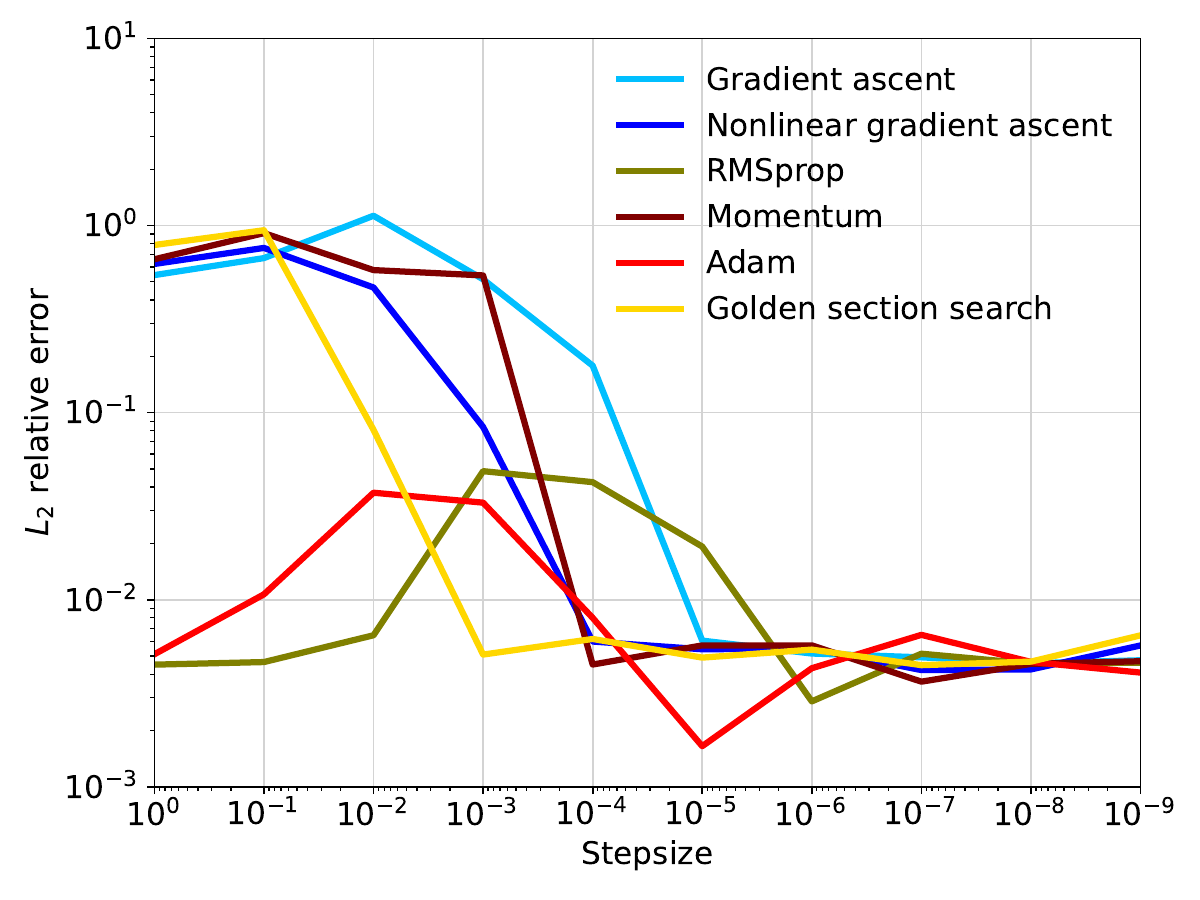}
         \caption{}
         \label{fig:Allen_Cahn_L2_vs_stepsize}
     \end{subfigure}
     \begin{subfigure}[b]{0.49\textwidth}
         \centering
         \includegraphics[width=\textwidth]{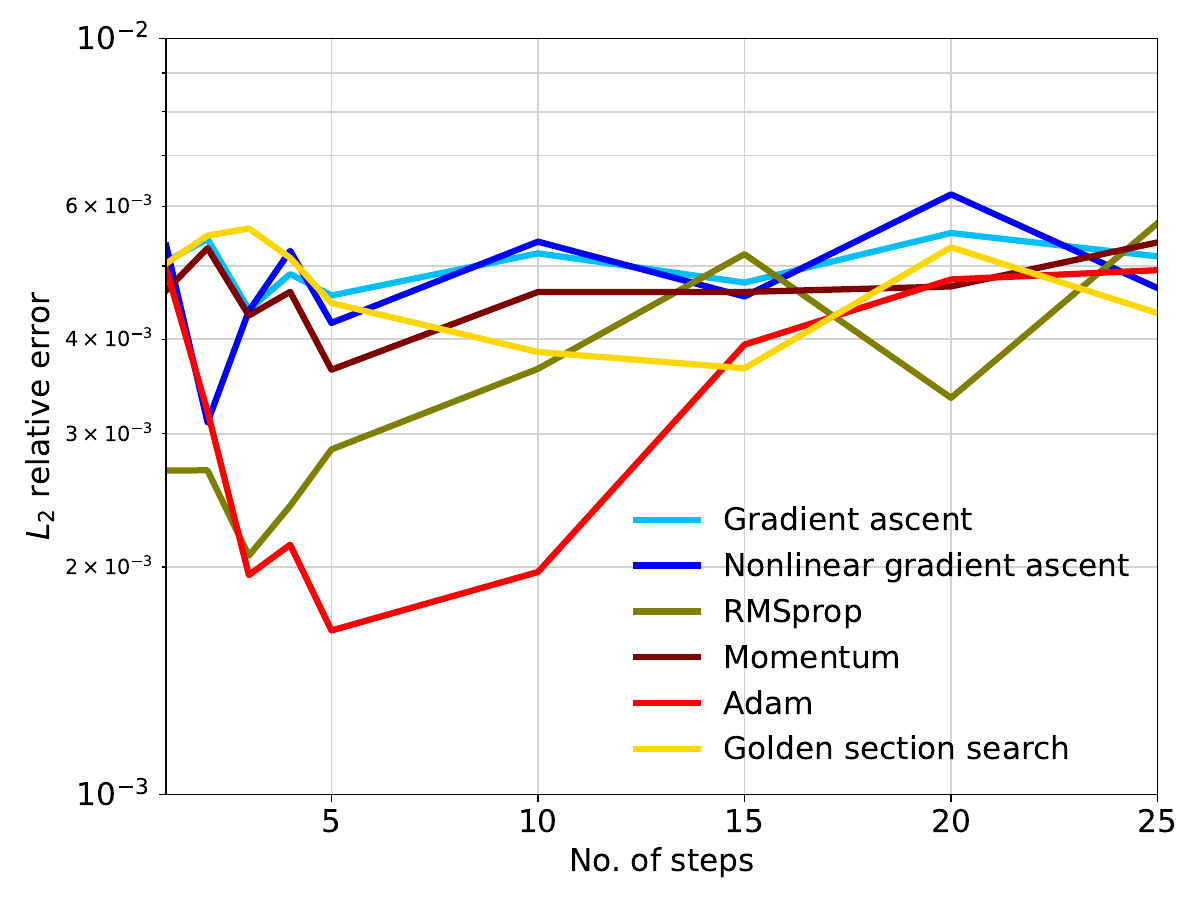}
         \caption{}
         \label{fig:Allen_Cahn_L2_vs_iter}
     \end{subfigure}
        \caption{Mean of the test error for PACMANN with the different optimization algorithms listed in~\Cref{sec:Methodology} for a varying (a) stepsize and (b) number of steps for the Allen-Cahn equation example.}
        \label{fig:Allen_Cahn_L2_vs_step_iter}
\end{figure}

In the following example, we consider the one-dimensional Allen-Cahn equation:
\begin{equation}
\label{eq:Allen_Cahn_equation_problem_statement}
\left\{
\begin{aligned}
u_t & = d u_{xx} + 5(u-u^3), \quad x \in [-1,1], \quad t \in [0,1], \\
u(x,0) & = x^2\cos(\pi x), \\
u(-1,t) & = u(1,t) = -1. \\
\end{aligned} \right.
\end{equation}
We choose a diffusion coefficient of $d = 0.001$. Similar to the Burgers' equation problem in the previous section, the number of collocation points is set to 2500, the number of boundary points to 80, and the number of initial condition points to 160. The network architecture used for this example consists of four hidden layers of 64 neurons. To compute the $L_2$ relative error, we compare the network prediction to a numerical solution generated using a spectral solver and made available via the DeepXDE library~\cite{Lu2021DeepXDE:Equations}.

The mean and standard deviation of the test error and the mean runtime for each sampling method are given in~\Cref{tab:Allen_Cahn_overview_L2_error_resamplers}. As for the Burgers' equation example, the static uniform grid fails to learn the solution satisfactorily, and the static Hammersley grid and random resampling offer a significant improvement in accuracy. Out of all sampling methods considered, PACMANN in combination with the Adam optimizer and a stepsize of $10^{-5}$ achieves the lowest test error. The next-best method, RAR-D, attains a lower prediction accuracy at the same computational cost. \Cref{tab:Allen_Cahn_overview_L2_error_optimizers} demonstrates that Adam results in the lowest mean and standard deviation of the test error for our method in comparison to the other optimization algorithms considered.

\Cref{fig:Allen_Cahn_L2_vs_collocation} shows the behavior of the different sampling methods when varying the number of collocation points. We find that PACMANN with the Adam optimizer converges the fastest and outperforms the other sampling methods for nearly all numbers of collocation points considered. Most sampling methods stagnate above 2500 collocation points, except for the static uniform grid, which requires up to 10\,000 collocation points to compete with the accuracy of the other sampling methods.

Moreover, \Cref{fig:Allen_Cahn_L2_vs_period} depicts the accuracy of the sampling methods considered for a range of resampling periods. In contrast to our findings in the previous Burgers' equation example, the resampling period has a reduced influence on the test error in this example. Notably, most sampling methods plateau with an increasing resampling period.

Furthermore, we test the behavior of PACMANN with the different optimization algorithms listed in~\Cref{sec:Methodology} by varying the stepsize and the number of steps hyperparameters. Similarly to the results in~\cref{sec:Burgers_equation} for the Burgers' equation, the stepsize hyperparameter has a significant influence on the test error, as demonstrated by~\Cref{fig:Allen_Cahn_L2_vs_stepsize}. We note, again, the formation of two groups, which occurs due to the random resampling of points that have been moved outside the domain. When using RMSprop and Adam combined with large stepsizes, all 2500 collocation points leave the domain. As a result, the test error converges toward the error of $0.42 \% \pm 0.28 \%$ found when using random resampling; cf.~\cref{tab:Allen_Cahn_overview_L2_error_resamplers}.

Finally, we also vary the number of optimization steps for resampling. \Cref{fig:Allen_Cahn_L2_vs_iter} illustrates that the test error of most optimization algorithms is not significantly impacted by changing the number of steps. \Cref{fig:Allen_Cahn_L2_vs_stepsize,fig:Allen_Cahn_L2_vs_iter} support our earlier finding that the stepsize plays a dominant role in the test error obtained compared to the limited influence of the number of steps.

Based on the observations made in the Burgers' and Allen-Cahn equation examples, we recommend fixing the number of steps to five, which removes a hyperparameter from PACMANN. Furthermore, we recommend five steps as it takes advantage of the increased accuracy achieved by our method with the Adam optimizer at multiple steps (see~\Cref{fig:Burgers_L2_vs_iter,fig:Allen_Cahn_L2_vs_iter}), while keeping the computational cost low; see~\Cref{fig:Burgers_L2_Duration_vs_iter}.

\subsection{2D Navier-Stokes equation (inverse)}

\begin{table}[t]
\centering
\begin{tabular}{@{}lccccc@{}}
\toprule
\multicolumn{1}{c}{\multirow{3}{*}{Sampling method}} & \multicolumn{4}{c}{$L_2$ relative error} & \multirow{3}{*}{\begin{tabular}[c]{@{}c@{}}Mean \\ runtime {[}s{]}\end{tabular}} \\ \cmidrule(lr){2-5}
\multicolumn{1}{c}{} & \multicolumn{2}{c}{$\lambda_1$} & \multicolumn{2}{c}{$\lambda_2$} &  \\ \cmidrule(lr){2-5}
\multicolumn{1}{c}{} & Mean & 1 SD & Mean & 1 SD &  \\ \midrule
\multicolumn{1}{l}{Uniform grid} & 0.05\% & \textbf{0.01\%} & 0.72\% & 0.43\% & 1506 \\
\multicolumn{1}{l}{Hammersley grid} & 0.08\% & 0.04\% & 0.89\% & 0.52\% & \textbf{1492} \\
\multicolumn{1}{l}{Random resampling} & 0.12\% & 0.05\% & 0.65\% & 0.46\% & 1514 \\ \cmidrule(r){1-6}
\multicolumn{1}{l}{RAR} & 0.30\% & 0.06\% & 1.44\% & 0.90\% & 1520 \\
\multicolumn{1}{l}{RAD} & 0.23\% & 0.06\% & 1.38\% & 0.79\% & 1583 \\
\multicolumn{1}{l}{RAR-D} & 0.08\% & 0.05\% & 0.84\% & 0.57\% & 1525 \\ \cmidrule(r){1-6}
\multicolumn{1}{l}{PACMANN-Adam} & \textbf{0.03\%} & 0.03\% & \textbf{0.53\%} & \textbf{0.19\%} & 1559 \\ \bottomrule
\end{tabular}
\caption{Overview of the mean and standard deviation of the test error for $\lambda_1$ and $\lambda_2$ and the mean runtime for each sampling method for the inverse Navier-Stokes equation example. The best result in each column is marked in boldface.}
\label{tab:Navier_Stokes_overview_L2_error}
\end{table}

Next, we consider an inverse problem based on the two-dimensional Navier-Stokes equation describing the flow of an incompressible fluid past a cylinder discussed by Raissi et al.~in~\cite{Raissi2019Physics-informedEquations}, given by:
\begin{equation}
\label{eq:Navier_Stokes_equation_problem_statement}
\left\{
\begin{aligned}
\frac{\partial \boldsymbol{v}}{\partial t} + \lambda_1 \boldsymbol{v} \cdot \nabla \boldsymbol{v} &= -\nabla p + \lambda_2 \nabla^{2} \boldsymbol{v}, \quad (x,y) \in [1,8] \times [-2,2], \quad t \in [0,7], \\
\nabla \cdot \boldsymbol{v} &= 0.
\end{aligned}
\right.
\end{equation}
Here, $u$ and $v$ are the $x$- and $y$-components of the velocity field, and $p$ denotes the pressure. The scalar parameter $\lambda_1$ scales the convective term, and $\lambda_2$ represents the dynamic (shear) viscosity. In this example, we are interested in learning the values of $\lambda_1$ and $\lambda_2$ based on a data set created by Raissi et al.~\cite{Raissi2019Physics-informedEquations} using a spectral solver. The data set contains the values of $u$, $v$, and $p$ determined for a large set of points $(x,y,t)$.
The true values of $\lambda_1$ and $\lambda_2$ are 1 and 0.01, respectively. For this inverse problem, we train the PINN on 7000 randomly selected points from this data set. In addition, we sample 700 collocation points, 200 points on the boundary condition, and 100 points on the initial condition. The network architecture consists of six hidden layers containing 50 neurons each.

The mean and standard deviation of the test error for both $\lambda_1$ and $\lambda_2$ and the mean runtime for each of the sampling methods are given in~\Cref{tab:Navier_Stokes_overview_L2_error}. In this inverse problem, PACMANN in combination with Adam at a stepsize of $10^{-2}$ achieves the lowest test error for $\lambda_1$ and $\lambda_2$ at a slightly higher computational cost compared to the second best adaptive method, RAR-D. Furthermore, we note that the non-adaptive sampling methods generally outperform the other adaptive methods considered in this study, both in the mean and standard deviation of the test error.

Qualitatively, we observe similar behavior of the different sampling methods when varying the number of collocation points and the resampling period for the 2D Navier-Stokes equation example and the following examples, as observed in the Burgers' equation example. This observation also applies to the behavior of PACMANN when changing the stepsize and the number of steps hyperparameters, namely, that the number of steps generally has a smaller impact on the prediction accuracy than the stepsize. Therefore, for the sake of conciseness, we do not repeat the analysis of the behavior of the different sampling methods for this two-dimensional Navier-Stokes equation example and the following examples.

\subsection{5D Poisson's equation}

\begin{table}[t]
\centering
\begin{tabular}{@{}lccccc@{}}
\toprule
\multicolumn{1}{c}{\multirow{2}{*}{Sampling method}} & \multicolumn{2}{c}{$L_2$ relative error} & \multicolumn{2}{c}{$H^1$ semi-norm} & \multirow{2}{*}{\begin{tabular}[c]{@{}c@{}}Mean\\ runtime {[}s{]}\end{tabular}} \\ \cmidrule(lr){2-5}
\multicolumn{1}{c}{}                                 & Mean               & 1 SD               & Mean              & 1 SD            &                                                                                 \\ \midrule
Uniform grid                                         & 40.32\%            & 1.18\%             & 62.58\%           & 1.74\%          & \textbf{744}                                                                    \\
Hammersley grid                                      & 82.64\%            & 2.95\%             & 90.48\%           & 2.20\%          & 752                                                                             \\
Random resampling                                    & 11.47\%            & 1.13\%             & 13.36\%           & 1.66\%          & 751                                                                             \\ \midrule
RAR                                                  & 62.00\%            & 7.13\%             & 69.97\%           & 6.62\%          & 773                                                                             \\
RAD                                                  & 11.69\%            & 1.78\%             & 13.19\%           & 1.31\%          & 841                                                                             \\
RAR-D                                                & 89.31\%            & 1.98\%             & 93.83\%           & 2.05\%          & 784                                                                             \\ \midrule
PACMANN-Adam                                                & \textbf{8.35\%}    & \textbf{0.54\%}    & \textbf{10.32\%}  & \textbf{0.43\%} & 786                                                                             \\ \bottomrule
\end{tabular}
\caption{Overview of the mean and standard deviation of the $L_2$ relative error and the $H^1$ semi-norm, and the mean runtime for each sampling method for the Poisson's equation example. The best result in each column is marked in boldface.}
\label{tab:Poisson_overview_L2_error}
\end{table}

In the following problem, we apply PACMANN to the Poisson equation in five dimensions:
\begin{equation}
\label{eq:Poissons_equation_problem_statement}
\left\{
\begin{aligned}
-\Delta v \left( \boldsymbol{x} \right) & = f \left( \boldsymbol{x} \right), & \quad & \boldsymbol{x} \in [-1,1]^5 \\
v \left( \boldsymbol{x} \right) & = 0, & & \boldsymbol{x} \in \partial \boldsymbol{\Omega}.\\
\end{aligned} \right.
\end{equation}
We choose the right-hand side function $f$ based on the manufactured solution
\begin{equation}
    \label{eq:Solution_Poissons_equation}
    v \left( \boldsymbol{x} \right) = \prod_{i=1}^5 \sin(\pi x_i),
\end{equation}
where $x_i$ is the $i$-th component of $\boldsymbol{x}$. For this example, we sample 750 collocation points and 750 points for the boundary condition. Moreover, the network architecture consists of four hidden layers of 64 neurons each.

The mean and standard deviation of the $L_2$ relative error and the $H^1$ semi-norm, and the mean runtime for each of the sampling methods are given in~\Cref{tab:Poisson_overview_L2_error}. In contrast to the previous examples, the adaptive methods RAR and RAR-D fail to improve the prediction accuracy. PACMANN in combination with the Adam optimizer and a stepsize of $10^{-2}$ achieves the lowest mean and standard deviation of the $L_2$ relative error and the $H^1$ semi-norm. Moreover, we point out the ability of this method to efficiently scale to high-dimensional problems.
We find that our method with the Adam optimizer is cheaper at a mean runtime of 786 seconds compared to RAD, the next-best adaptive sampling method, at 841 seconds.

\subsection{3D Navier-Stokes equation}

\begin{table}[htbp]
\centering
\begin{tabular}{@{}lccccccc@{}}
\toprule
\multicolumn{1}{c}{\multirow{3}{*}{Sampling method}} & \multicolumn{6}{c}{$L_2$ relative error}                                                                                                                             & \multicolumn{1}{c}{\multirow{3}{*}{\begin{tabular}[c]{@{}c@{}}Mean\\ runtime {[}s{]}\end{tabular}}} \\ \cmidrule(lr){2-7}
\multicolumn{1}{c}{}                                 & \multicolumn{2}{c}{$u$}                             & \multicolumn{2}{c}{$v$}                             & \multicolumn{2}{c}{$w$}                             & \multicolumn{1}{c}{}                                                                                \\ \cmidrule(lr){2-7}
\multicolumn{1}{c}{}                                 & \multicolumn{1}{c}{Mean} & \multicolumn{1}{c}{1 SD} & \multicolumn{1}{c}{Mean} & \multicolumn{1}{c}{1 SD} & \multicolumn{1}{c}{Mean} & \multicolumn{1}{c}{1 SD} & \multicolumn{1}{c}{}                                                                                \\ \midrule
Uniform grid                                         & 1.01\%                   & 0.16\%                   & 1.20\%                   & 0.34\%                   & 1.27\%                   & 0.41\%                   & 1506                                                                                                \\
Hammersley grid                                      & 0.67\%                   & 0.10\%                   & 0.68\%                   & 0.09\%                   & 0.48\%                   & 0.10\%                   & 1503                                                                                                \\
Random resampling                                    & 0.50\%                   & 0.09\%                   & \textbf{0.49\%}          & 0.09\%                   & 0.37\%                   & 0.07\%                   & 1467                                                                                                \\ \midrule
RAR                                                  & 0.59\%                   & 0.13\%                   & 0.56\%                   & 0.07\%                   & 0.42\%                   & 0.05\%                   & \textbf{1462}                                                                                       \\
RAD                                                  & 0.49\%                   & 0.06\%                   & 0.51\%                   & \textbf{0.06\%}          & 0.37\%                   & \textbf{0.03\%}          & 1531                                                                                                \\
RAR-D                                                & 0.65\%                   & 0.12\%                   & 0.68\%                   & 0.10\%                   & 0.47\%                   & 0.06\%                   & 1497                                                                                                \\ \midrule
PACMANN-Adam                                                & \textbf{0.46\%}          & \textbf{0.04\%}          & \textbf{0.49\%}          & \textbf{0.09\%}          & \textbf{0.33\%}          & 0.05\%                   & 1524                                                                                                \\ \bottomrule
\end{tabular}
\caption{Overview of the mean and standard deviation of the $L_2$ relative error for $u$, $v$, and $w$ and the mean runtime for each sampling method for the 3D Navier-Stokes equation example. The best result in each column is marked in boldface.}
\label{tab:3D_Navier_Stokes_overview_L2_error}
\end{table}

\begin{table}[htbp]
\centering
\begin{tabular}{@{}lccccccc@{}}
\toprule
\multicolumn{1}{c}{\multirow{3}{*}{Sampling method}} & \multicolumn{6}{c}{$H^1$ semi-norm}                                                                                                                             & \multicolumn{1}{c}{\multirow{3}{*}{\begin{tabular}[c]{@{}c@{}}Mean\\ runtime {[}s{]}\end{tabular}}} \\ \cmidrule(lr){2-7}
\multicolumn{1}{c}{}                                 & \multicolumn{2}{c}{$u$}                             & \multicolumn{2}{c}{$v$}                             & \multicolumn{2}{c}{$w$}                             & \multicolumn{1}{c}{}                                                                                \\ \cmidrule(lr){2-7}
\multicolumn{1}{c}{}                                 & \multicolumn{1}{c}{Mean} & \multicolumn{1}{c}{1 SD} & \multicolumn{1}{c}{Mean} & \multicolumn{1}{c}{1 SD} & \multicolumn{1}{c}{Mean} & \multicolumn{1}{c}{1 SD} & \multicolumn{1}{c}{}                                                                                \\ \midrule
Uniform grid                                         & 1.30\%                   & 0.18\%                   & 1.40\%                   & 0.19\%                   & 1.24\%                   & 0.26\%                   & 1506                                                                                                \\
Hammersley grid                                      & 0.43\%                   & 0.05\%                   & 0.45\%                   & 0.08\%                   & 0.29\%                   & 0.05\%                   & 1503                                                                                                \\
Random resampling                                    & 0.35\%                   & 0.04\%                   & \textbf{0.35\%}          & 0.05\%                   & 0.25\%                   & 0.03\%                   & 1467                                                                                                \\ \midrule
RAR                                                  & 0.43\%                   & 0.06\%                   & 0.42\%                   & \textbf{0.04\%}                   & 0.29\%                   & 0.03\%                   & \textbf{1462}                                                                                       \\
RAD                                                  & 0.36\%                   & 0.04\%                   & 0.38\%                   & 0.06\%                   & 0.24\%                   & \textbf{0.02\%}          & 1531                                                                                                \\
RAR-D                                                & 0.43\%                   & 0.07\%                   & 0.42\%                   & 0.05\%                   & 0.27\%                   & 0.04\%                   & 1497                                                                                                \\ \midrule
PACMANN-Adam                                                & \textbf{0.33\%}          & \textbf{0.03\%}          & \textbf{0.35\%}          & \textbf{0.04\%}          & \textbf{0.22\%}          & 0.04\%                   & 1524                                                                                                \\ \bottomrule
\end{tabular}
\caption{Overview of the mean and standard deviation of the $H^1$ semi-norm for $u$, $v$, and $w$ and the mean runtime for each sampling method for the 3D Navier-Stokes equation example. The best result in each column is marked in boldface.}
\label{tab:3D_Navier_Stokes_overview_H1_seminorm}
\end{table}

We next apply PACMANN to the incompressible Navier-Stokes equation in a cube:
\begin{equation}
\label{eq:Navier_Stokes_equation_cube_problem_statement}
\left\{
\begin{aligned}
\frac{\partial \boldsymbol{v}}{\partial t} + Re \: \boldsymbol{v} \cdot \nabla \boldsymbol{v} +Re \: \nabla p - \nabla^{2} \boldsymbol{v} &= f(\boldsymbol{x}), \quad (x,y,z) \in [-1,1]^3, \quad t \in [0,1], \\
\nabla \cdot \boldsymbol{v} &= 0.
\end{aligned}
\right.
\end{equation}
Here, we set the Reynolds number to $Re = 1$. Furthermore, we choose the right-hand side function $f$ based on the manufactured solution
\begin{equation}
\label{eq:Navier_Stokes_equation_cube_manufactured_solution}
\begin{aligned}
u &= \cos(x) \sin(y) \sin(z) e^{-t}\\
v &= \sin(x) \cos(y) \sin(z) e^{-t}\\
w &= -2 \cos(x) \cos(y) \cos(z) e^{-t}.\\
\end{aligned}
\end{equation}
We compare PACMANN against the other sampling methods for the Navier-Stokes problem with Dirichlet boundary conditions and initial conditions prescribed by the manufactured solution. For this problem, we sample 300 collocation points, and we employ 75 points for the boundary conditions and 75 points for the initial conditions. The network architecture used for this problem consists of four hidden layers of 64 neurons.

The mean and standard deviation of the $L_2$ relative error and the $H^1$ semi-norm for $u, v$, and $w$, and the mean runtime for each of the sampling methods are given in~\Cref{tab:3D_Navier_Stokes_overview_L2_error,tab:3D_Navier_Stokes_overview_H1_seminorm}. In this problem, PACMANN in combination with Adam at a stepsize of $10^{-4}$ achieves the lowest errors for $u$, $v$, and $w$ while maintaining a computational cost comparable to that of other adaptive sampling methods.

\subsection{Re-entrant corner in a disk}

\begin{figure}[t]
    \centering
    \begin{subfigure}{0.32\textwidth}
        \centering
        \includegraphics[width=0.96\linewidth,trim=7mm 18mm 9mm 18mm,clip]{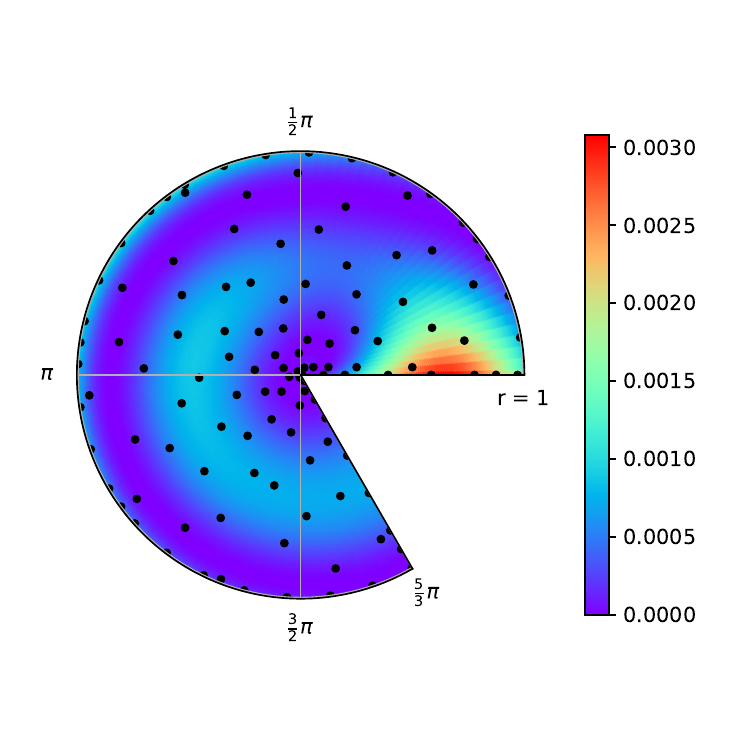}
        \caption{}
        \label{fig:sub1}
    \end{subfigure}
    \begin{subfigure}{0.32\textwidth}
        \centering
        \includegraphics[width=0.96\linewidth,trim=7mm 18mm 9mm 18mm,clip]{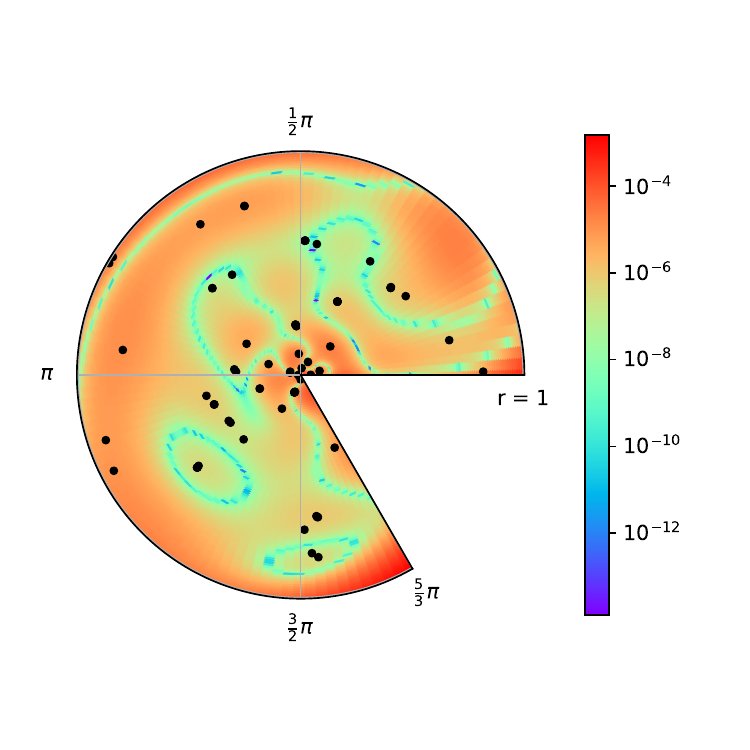}
        \caption{}
        \label{fig:sub2}
    \end{subfigure}
    \begin{subfigure}{0.32\textwidth}
        \centering
        \includegraphics[width=0.96\linewidth,trim=7mm 18mm 9mm 18mm,clip]{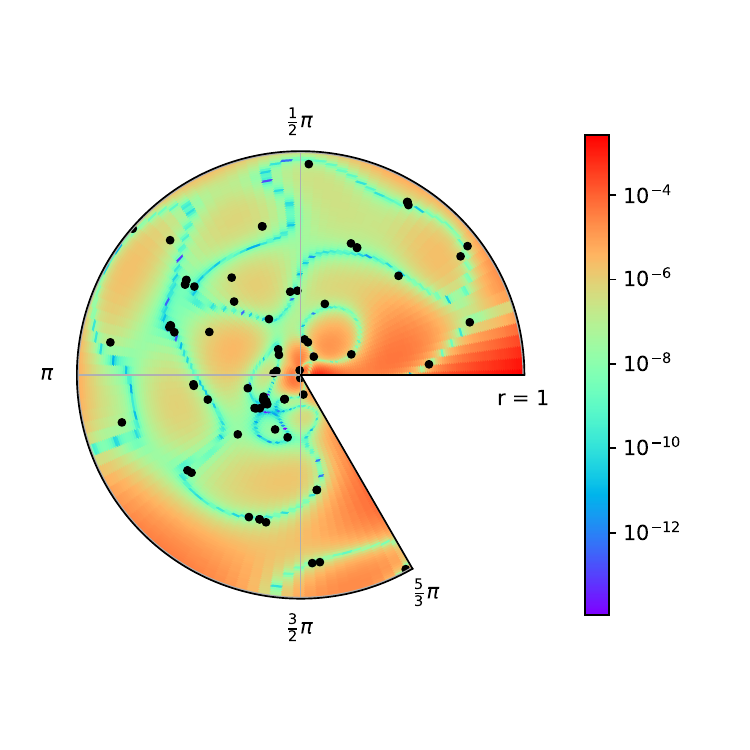}
        \caption{}
        \label{fig:sub3}
    \end{subfigure}

    \caption{Location of the collocation points (a) before training, (b) mid-training at 25\,000 iterations, and (c) after training using PACMANN with Adam for the Laplace equation over a re-entrant corner in a disk example. The color indicates the values of the squared residual.}
    \label{fig:three_residual_plots}
\end{figure}

\begin{table}[t]
\centering
\begin{tabular}{@{}lccccc@{}}
\toprule
\multicolumn{1}{c}{\multirow{2}{*}{Sampling method}} & \multicolumn{2}{c}{$L_2$ relative error} & \multicolumn{2}{c}{$H^1$ semi-norm} & \multirow{2}{*}{\begin{tabular}[c]{@{}c@{}}Mean\\ runtime {[}s{]}\end{tabular}} \\ \cmidrule(lr){2-5}
\multicolumn{1}{c}{}                                 & Mean               & 1 SD               & Mean             & 1 SD             &                                                                                 \\ \midrule
Uniform grid                                         & 6.32\%             & 1.46\%             & 32.36\%         & 6.50\%           & \textbf{514}                                                                    \\
Hammersley grid                                      & 0.47\%             & \textbf{0.17\%}    & 5.90\%           & 0.91\%           & 539                                                                             \\
Random resampling                                    & 0.49\%             & 0.18\%             & 5.45\%           & \textbf{0.82\%}  & 536                                                                             \\ \midrule
RAR                                                  & 0.58\%             & \textbf{0.17\%}    & 6.12\%           & 1.00\%           & 538                                                                             \\
RAD                                                  & 0.48\%             & 0.25\%             & \textbf{4.46\%}  & 1.31\%           & 594                                                                             \\
RAR-D                                                & 0.96\%             & 0.74\%             & 9.39\%           & 6.27\%           & 543                                                                             \\ \midrule
PACMANN-Adam                                                & \textbf{0.43\%}    & \textbf{0.17\%}    & 5.00\%           & 1.03\%           & 544                                                                             \\ \bottomrule
\end{tabular}
\caption{Overview of the mean and standard deviation of the $L_2$ relative error and the $H^1$ semi-norm, and the mean runtime for each sampling method for the Laplace equation over a re-entrant corner in a disk example. The best result in each column is marked in boldface.}
\label{tab:2D_Laplace_equation_overview_L2_error}
\end{table}

Finally, we consider a problem involving a domain with a re-entrant corner, where the solution is expected to be less regular. Specifically, we consider the Laplace equation in polar coordinates over a sector:
\begin{equation}
\label{eq:Laplace_disk_problem_statement}
\left\{
\begin{aligned}
u_{rr} + \frac{1}{r}u_{r} + \frac{1}{r^2} u_{\theta \theta} &= 0, & \quad & r \in [0, 1], \quad \theta \in \left[0, \frac{5\pi}{3}\right], \\
u(r, \theta) & = h(r, \theta), & & (r,\theta) \in \partial \boldsymbol{\Omega}.
\end{aligned}
\right.
\end{equation}
The exact solution is used for the boundary condition $h(r,\theta)$ and is given by
\begin{equation*}
    h(r,\theta) = r^{\frac{3}{5}}\sin\left(\frac{3}{5}\theta\right),
\end{equation*}
which has an algebraic singularity at the origin; cf.~\cite[pp.~110--116]{Dauge1988EllipticDomains}. In this example, we sample 75 collocation points and we set the number of points for the boundary condition to 75. The network architecture used for this problem consists of four hidden layers of 64 neurons.

The mean and standard deviation of the $L_2$ relative error and the $H^1$ semi-norm, and the mean runtime for each of the sampling methods are given in~\Cref{tab:2D_Laplace_equation_overview_L2_error}. In this example, PACMANN in combination with the Adam optimizer at a stepsize of $10^{-2}$ achieves the lowest $L_2$ relative error at a computational cost comparable to that of the adaptive sampling methods RAR and RAR-D, while RAD achieves the lowest $H^1$ semi-norm. Furthermore, we observe a collective movement of points toward the origin, where the algebraic singularity is located; see~\Cref{fig:three_residual_plots}.

\section{Conclusions}
\label{sec:Conclusions}
In this work, we presented the Point Adaptive Collocation Method for Artificial Neural Networks (PACMANN), a novel adaptive collocation point sampling method for physics-informed neural networks. This approach uses the gradient of the physics-informed loss terms, that is, of the squared residual, as guiding information to move collocation points toward areas of large residuals. The problem of moving points is formulated as a maximization problem, which can be approached using an optimization algorithm of choice, such as gradient ascent or Adam. Points are moved several times while training is halted. Our approach can be tuned using three additional hyperparameters, namely the resampling period, the size of the step used to move the collocation points, and the number of times that points are moved while training is halted. While this work demonstrates PACMANN for PINNs, we note that the method can also be applied to other collocation-based approaches.

We studied the sensitivity of our method to these hyperparameters, and we observed that the stepsize has a particularly large impact on the solution accuracy. Conversely, we found that five iteration steps are sufficient to achieve a good balance between accuracy and efficiency. We then investigated the accuracy and efficiency of PACMANN in combination with various optimization algorithms and concluded that the Adam optimizer performs the best.

Furthermore, we compared the performance of PACMANN to existing state-of-the-art adaptive and non-adaptive collocation approaches, including random resampling, RAR, and RAD, and demonstrated that our method achieves state-of-the-art performance in terms of the accuracy/efficiency tradeoff for lower-dimensional benchmarks, while outperforming the state-of-the-art for high-dimensional problems and the case of an irregular domain. In addition, we showed the effectiveness of our approach in solving inverse problems.

In particular, the results of the numerical experiments demonstrate that PACMANN achieves high prediction accuracy and efficiently scales to higher dimensions without introducing significant computational overhead.

\section*{Acknowledgments}
The authors acknowledge the use of computational resources of the DelftBlue supercomputer, provided by Delft High Performance Computing Centre (\url{https://www.tudelft.nl/dhpc}).

\section*{Declaration of generative AI and AI-assisted technologies in the writing process}
During the preparation of this work the authors used ChatGPT (OpenAI) in order to improve the readability and language of the manuscript. After using this tool, the authors reviewed and edited the content as needed and take full responsibility for the content of the published article.

\appendix
\section{Relative error behavior of all optimization algorithms considered}
\label{sec:Appendix_A}

\begin{figure}[t]
     \centering
     \begin{subfigure}[b]{0.49\textwidth}
         \centering
         \includegraphics[width=\textwidth]{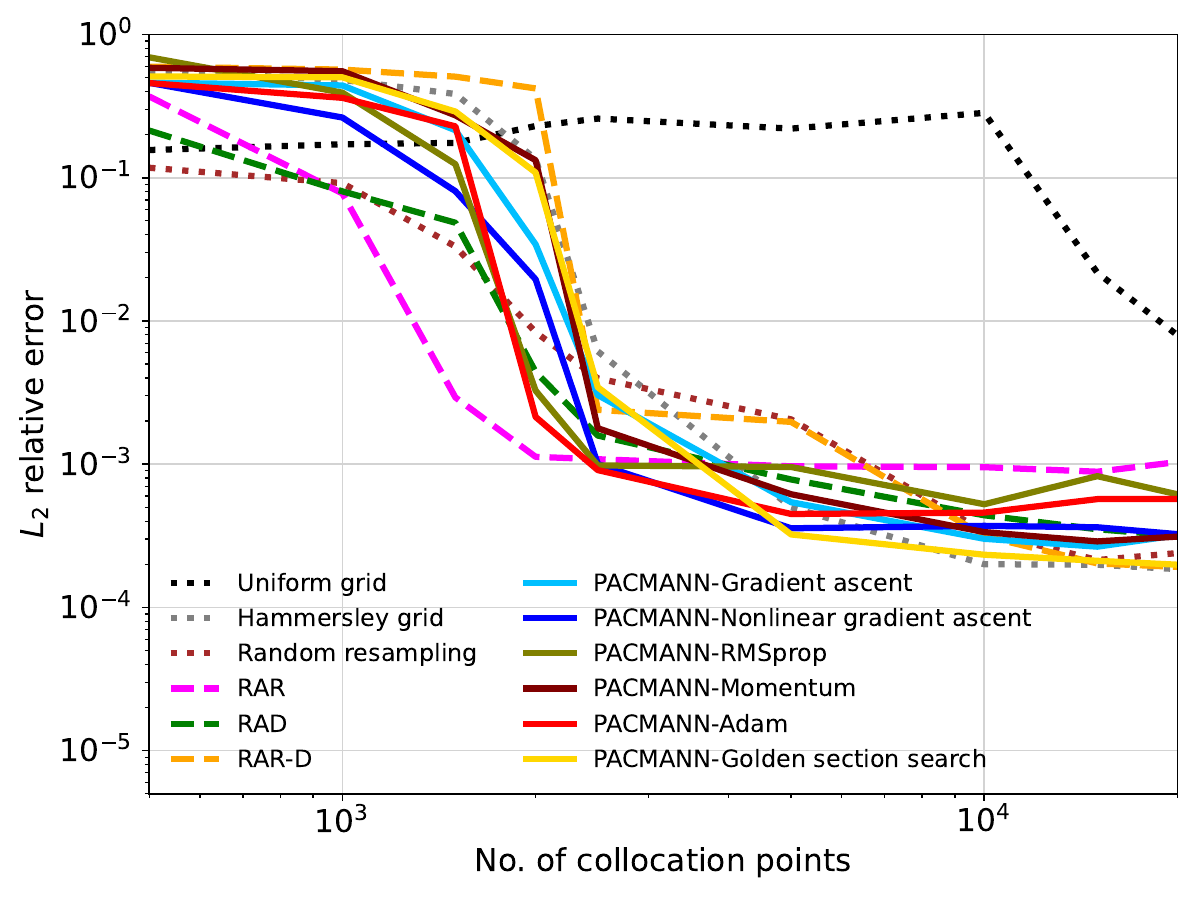}
         \caption{}
         \label{fig:Appendix_Burgers_L2_vs_collocation}
     \end{subfigure}
     \begin{subfigure}[b]{0.49\textwidth}
         \centering
         \includegraphics[width=\textwidth]{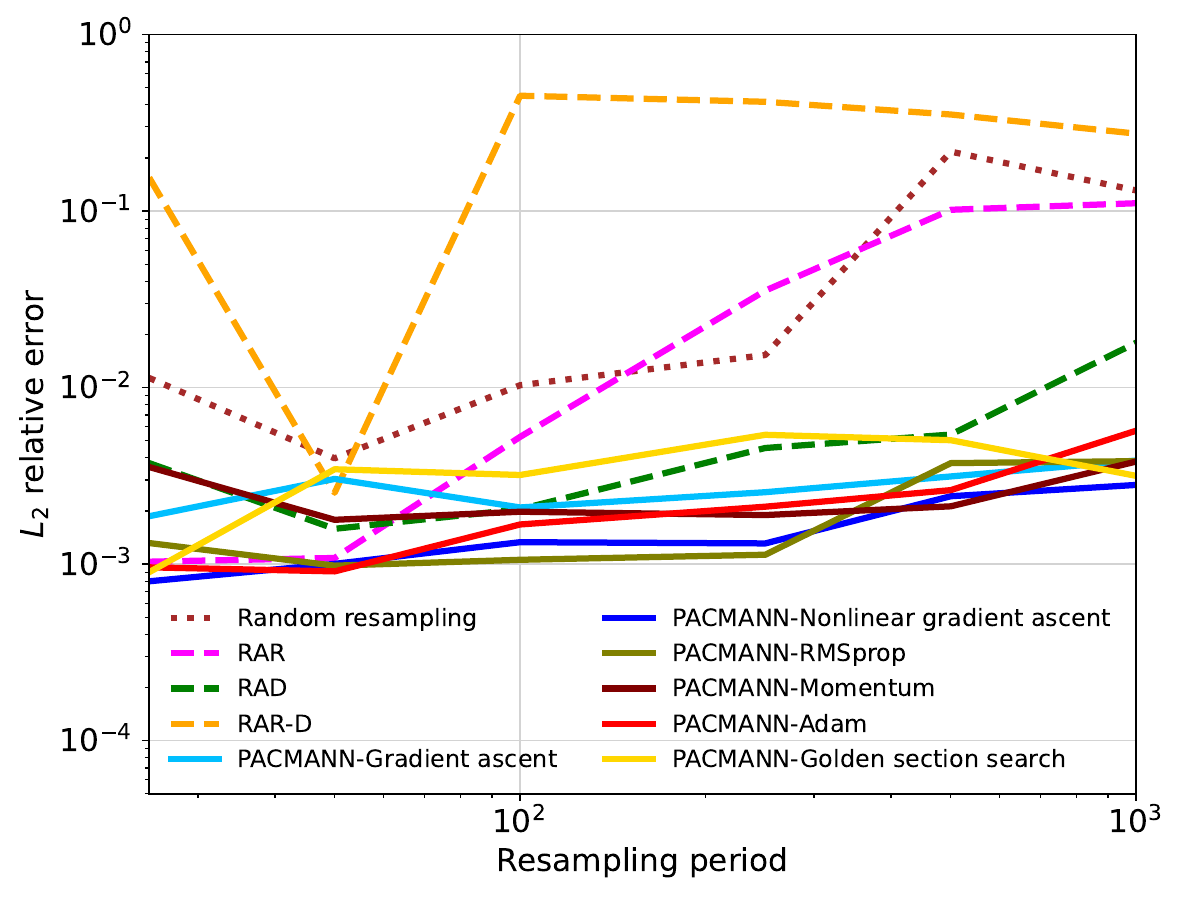}
         \caption{}
         \label{fig:Appendix_Burgers_L2_vs_period}
     \end{subfigure}
        \caption{Mean of the test error for each of the sampling methods and optimization algorithms considered for a varying (a) number of collocation points and (b) resampling period for the Burgers' equation example described in~\Cref{sec:Results}.}
        \label{fig:Appendix_Burgers_L2_vs_main_hyperparam}
\end{figure}

This appendix contains figures which describe the behavior of PACMANN in combination with all optimization methods listed in~\Cref{sec:Methodology}. \Cref{fig:Appendix_Burgers_L2_vs_collocation} illustrates the prediction accuracy of various optimization algorithms across a range of 500 to 20\,000 collocation points for the Burgers' equation example in~\Cref{sec:Burgers_equation}. Overall, the majority of algorithms demonstrate equivalent performance in terms of accuracy. Notably, the golden section search algorithm achieves the lowest error at large numbers of collocation points (20\,000 points). Furthermore, we evaluate the behavior of the optimization algorithms for a varying resampling period, from 25 to 1000 iterations, as depicted in~\Cref{fig:Appendix_Burgers_L2_vs_period}. Most optimization algorithms exhibit a behavior similar to Adam, as highlighted in~\Cref{sec:Results}. Notably, PACMANN loses prediction accuracy slower for increasing resampling period compared to the other sampling methods.

\begin{figure}[t]
     \centering
     \begin{subfigure}[b]{0.49\textwidth}
         \centering
         \includegraphics[width=\textwidth]{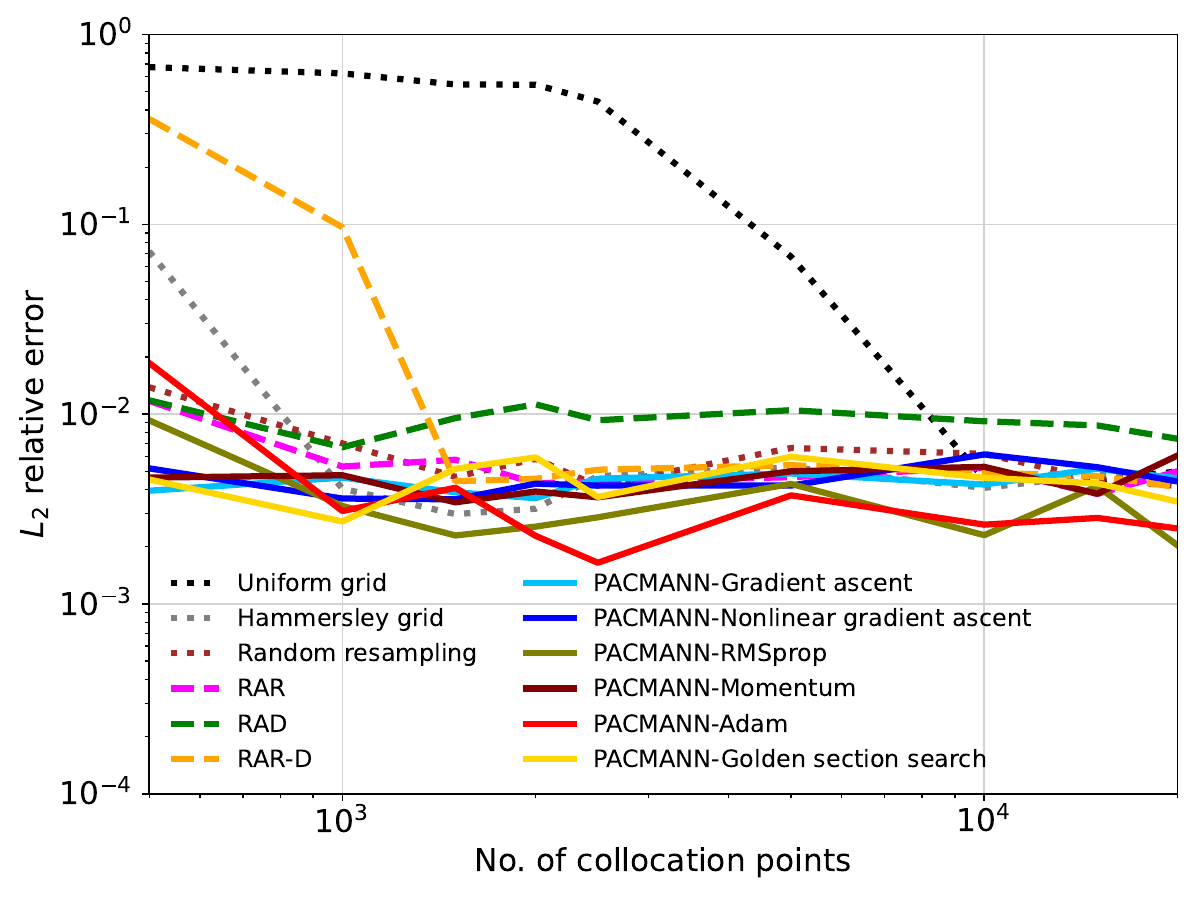}
         \caption{}
         \label{fig:Appendix_Allen_Cahn_L2_vs_collocation}
     \end{subfigure}
     \begin{subfigure}[b]{0.49\textwidth}
         \centering
         \includegraphics[width=\textwidth]{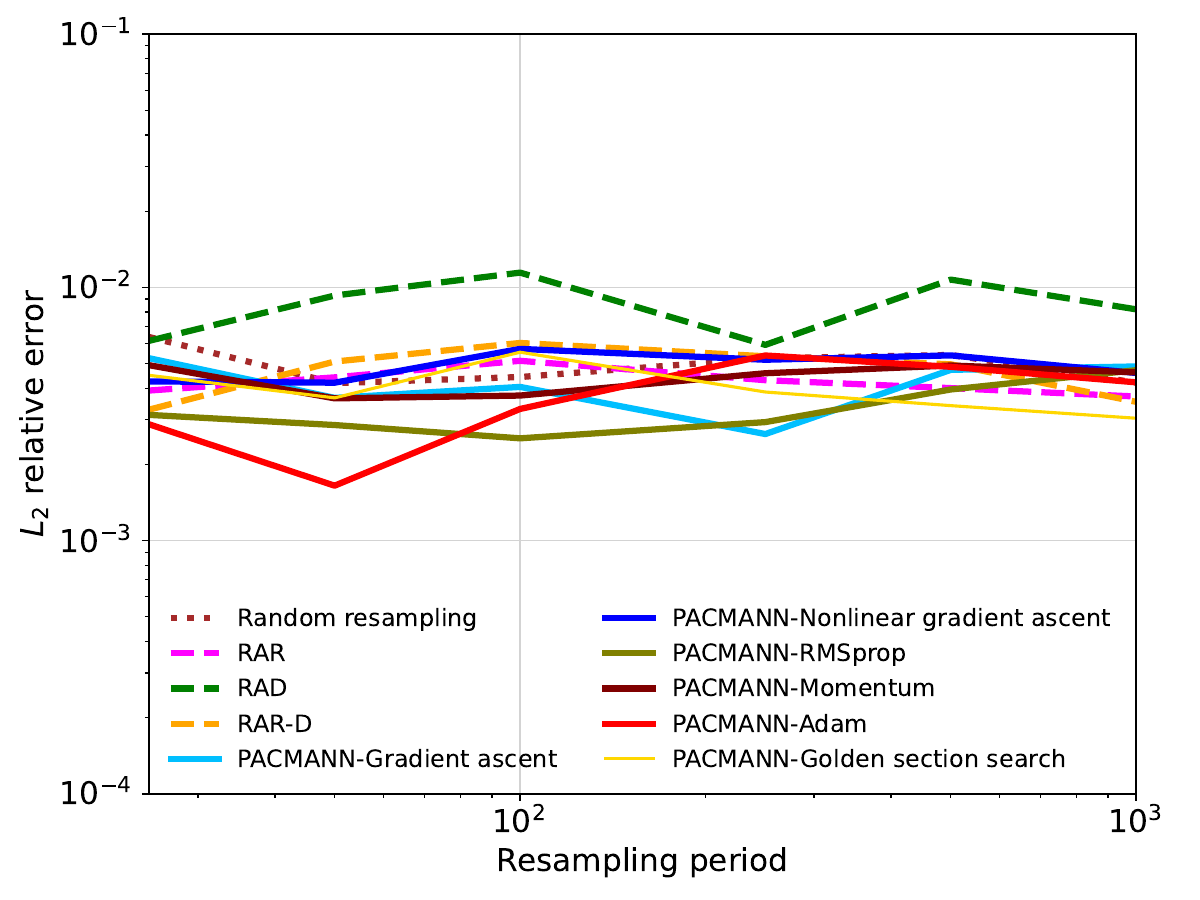}
         \caption{}
         \label{fig:Appendix_Allen_Cahn_L2_vs_period}
     \end{subfigure}
        \caption{Mean of the test error for each of the sampling methods and optimization algorithms considered for a varying (a) number of collocation points and (b) resampling period for the Allen-Cahn equation example described in~\Cref{sec:Results}.}
        \label{fig:Appendix_Allen_Cahn_L2_vs_main_hyperparam}
\end{figure}

Similarly, we vary the same hyperparameters for the Allen-Cahn equation example in~\cref{sec:Allen_Cahn_equation}. As shown in~\Cref{fig:Appendix_Allen_Cahn_L2_vs_collocation}, the other optimization algorithms display behavior comparable to Adam, with prediction accuracy remaining largely unaffected by an increasing number of collocation points. A similar trend is observed for the resampling period in~\Cref{fig:Appendix_Allen_Cahn_L2_vs_period}, where the accuracy shows minimal variation as the number of iterations increases.

\bibliographystyle{elsarticle-num}
\bibliography{Bibliography}

\begin{thebibliography}{10}
\expandafter\ifx\csname url\endcsname\relax
  \def\url#1{\texttt{#1}}\fi
\expandafter\ifx\csname urlprefix\endcsname\relax\def\urlprefix{URL }\fi
\expandafter\ifx\csname href\endcsname\relax
  \def\href#1#2{#2} \def\path#1{#1}\fi

\bibitem{Cybenko1989ApproximationFunction}
G.~Cybenko, {Approximation by superpositions of a sigmoidal function}, Mathematics of Control, Signals, and Systems 2~(4) (1989) 303--314.
\newblock \href {https://doi.org/10.1007/BF02551274} {\path{doi:10.1007/BF02551274}}.

\bibitem{Hornik1989MultilayerApproximators}
K.~Hornik, M.~Stinchcombe, H.~White, {Multilayer feedforward networks are universal approximators}, Neural Networks 2~(5) (1989) 359--366.
\newblock \href {https://doi.org/10.1016/0893-6080(89)90020-8} {\path{doi:10.1016/0893-6080(89)90020-8}}.

\bibitem{Dissanayake1994Neural-network-basedEquations}
M.~W. Dissanayake, N.~Phan‐Thien, {Neural-network-based approximations for solving partial differential equations}, Communications in Numerical Methods in Engineering 10~(3) (1994) 195--201.
\newblock \href {https://doi.org/10.1002/CNM.1640100303} {\path{doi:10.1002/CNM.1640100303}}.

\bibitem{Lagaris1998ArtificialEquations}
I.~E. Lagaris, A.~Likas, D.~I. Fotiadis, {Artificial neural networks for solving ordinary and partial differential equations}, IEEE Transactions on Neural Networks 9~(5) (1998) 987--1000.
\newblock \href {https://doi.org/10.1109/72.712178} {\path{doi:10.1109/72.712178}}.

\bibitem{Baydin2015AutomaticSurvey}
A.~G. Baydin, B.~A. Pearlmutter, A.~A. Radul, J.~M. Siskind, {Automatic differentiation in machine learning: a survey}, Journal of Machine Learning Research 18 (2015) 1--43.
\newblock \href {https://doi.org/10.5555/3122009.3242010} {\path{doi:10.5555/3122009.3242010}}.

\bibitem{Raissi2019Physics-informedEquations}
M.~Raissi, P.~Perdikaris, G.~E. Karniadakis, {Physics-informed neural networks: A deep learning framework for solving forward and inverse problems involving nonlinear partial differential equations}, Journal of Computational Physics 378 (2019) 686--707.
\newblock \href {https://doi.org/10.1016/J.JCP.2018.10.045} {\path{doi:10.1016/J.JCP.2018.10.045}}.

\bibitem{Cuomo2022ScientificNext}
S.~Cuomo, V.~S. Di~Cola, F.~Giampaolo, G.~Rozza, M.~Raissi, F.~Piccialli, {Scientific Machine Learning Through Physics–Informed Neural Networks: Where we are and What’s Next}, Journal of Scientific Computing 92~(3) (2022) 1--62.
\newblock \href {https://doi.org/10.1007/S10915-022-01939-Z} {\path{doi:10.1007/S10915-022-01939-Z}}.

\bibitem{Karniadakis2021Physics-informedLearning}
G.~E. Karniadakis, I.~G. Kevrekidis, L.~Lu, P.~Perdikaris, S.~Wang, L.~Yang, {Physics-informed machine learning}, Nature Reviews Physics 3~(6) (2021) 422--440.
\newblock \href {https://doi.org/10.1038/s42254-021-00314-5} {\path{doi:10.1038/s42254-021-00314-5}}.

\bibitem{Mao2020Physics-informedFlows}
Z.~Mao, A.~D. Jagtap, G.~E. Karniadakis, {Physics-informed neural networks for high-speed flows}, Computer Methods in Applied Mechanics and Engineering 360 (2020) 112789.
\newblock \href {https://doi.org/10.1016/J.CMA.2019.112789} {\path{doi:10.1016/J.CMA.2019.112789}}.

\bibitem{Cai2021FlowNetworks}
S.~Cai, Z.~Wang, F.~Fuest, Y.-J. Jeon, C.~Gray, G.~E. Karniadakis, {Flow over an espresso cup: Inferring 3D velocity and pressure fields from tomographic background oriented schlieren videos via physics-informed neural networks}, Journal of Fluid Mechanics 915 (2021).
\newblock \href {https://doi.org/10.1017/jfm.2021.135} {\path{doi:10.1017/jfm.2021.135}}.

\bibitem{Jin2021NSFnetsEquations}
X.~Jin, S.~Cai, H.~Li, G.~E. Karniadakis, {NSFnets (Navier-Stokes flow nets): Physics-informed neural networks for the incompressible Navier-Stokes equations}, Journal of Computational Physics 426 (2021) 109951.
\newblock \href {https://doi.org/10.1016/J.JCP.2020.109951} {\path{doi:10.1016/J.JCP.2020.109951}}.

\bibitem{Cai2021Physics-informedProblems}
S.~Cai, Z.~Wang, S.~Wang, P.~Perdikaris, G.~E. Karniadakis, {Physics-informed neural networks for heat transfer problems}, Journal of Heat Transfer 143~(6) (2021).
\newblock \href {https://doi.org/10.1115/1.4050542/1104439} {\path{doi:10.1115/1.4050542/1104439}}.

\bibitem{AminiNiaki2021Physics-informedManufacture}
S.~Amini~Niaki, E.~Haghighat, T.~Campbell, A.~Poursartip, R.~Vaziri, {Physics-informed neural network for modelling the thermochemical curing process of composite-tool systems during manufacture}, Computer Methods in Applied Mechanics and Engineering 384 (2021) 113959.
\newblock \href {https://doi.org/10.1016/J.CMA.2021.113959} {\path{doi:10.1016/J.CMA.2021.113959}}.

\bibitem{Shukla2020Physics-informedCracks}
K.~Shukla, P.~C. Di~Leoni, J.~Blackshire, D.~Sparkman, G.~E. Karniadakis, {Physics-informed neural network for ultrasound nondestructive quantification of surface breaking cracks}, Journal of Nondestructive Evaluation 39~(3) (2020).
\newblock \href {https://doi.org/10.1007/s10921-020-00705-1} {\path{doi:10.1007/s10921-020-00705-1}}.

\bibitem{Zhang2022AnalysesNetworks}
E.~Zhang, M.~Dao, G.~E. Karniadakis, S.~Suresh, {Analyses of internal structures and defects in materials using physics-informed neural networks}, Science Advances 8~(7) (2022).
\newblock \href {https://doi.org/10.1126/SCIADV.ABK0644} {\path{doi:10.1126/SCIADV.ABK0644}}.

\bibitem{Kovacs2022ConditionalNetworks}
A.~Kovacs, L.~Exl, A.~Kornell, J.~Fischbacher, M.~Hovorka, M.~Gusenbauer, L.~Breth, H.~Oezelt, M.~Yano, N.~Sakuma, A.~Kinoshita, T.~Shoji, A.~Kato, T.~Schrefl, {Conditional physics informed neural networks}, Communications in Nonlinear Science and Numerical Simulation 104 (2022) 106041.
\newblock \href {https://doi.org/10.1016/J.CNSNS.2021.106041} {\path{doi:10.1016/J.CNSNS.2021.106041}}.

\bibitem{Son2023AMotor}
S.~Son, H.~Lee, D.~Jeong, K.~Y. Oh, K.~Ho~Sun, {A novel physics-informed neural network for modeling electromagnetism of a permanent magnet synchronous motor}, Advanced Engineering Informatics 57 (2023) 102035.
\newblock \href {https://doi.org/10.1016/J.AEI.2023.102035} {\path{doi:10.1016/J.AEI.2023.102035}}.

\bibitem{Lu2021DeepXDE:Equations}
L.~Lu, X.~Meng, Z.~Mao, G.~E. Karniadakis, {DeepXDE: A deep learning library for solving differential equations}, SIAM Review 63~(1) (2021) 208--228.
\newblock \href {https://doi.org/10.1137/19M1274067} {\path{doi:10.1137/19M1274067}}.

\bibitem{Nabian2021EfficientSampling}
M.~A. Nabian, R.~J. Gladstone, H.~Meidani, {Efficient training of physics-informed neural networks via importance sampling}, Computer-Aided Civil and Infrastructure Engineering 36~(8) (2021) 962--977.
\newblock \href {https://doi.org/10.1111/mice.12685} {\path{doi:10.1111/mice.12685}}.

\bibitem{Wu2023ANetworks}
C.~Wu, M.~Zhu, Q.~Tan, Y.~Kartha, L.~Lu, {A comprehensive study of non-adaptive and residual-based adaptive sampling for physics-informed neural networks}, Computer Methods in Applied Mechanics and Engineering 403 (2023) 115671.
\newblock \href {https://doi.org/10.1016/J.CMA.2022.115671} {\path{doi:10.1016/J.CMA.2022.115671}}.

\bibitem{Guo2024TCAS-PINN:Method}
J.~Guo, H.~Wang, S.~Gu, C.~Hou, {TCAS-PINN: Physics-informed neural networks with a novel temporal causality-based adaptive sampling method}, Chinese Physics B 33~(5) (2024) 050701.
\newblock \href {https://doi.org/10.1088/1674-1056/AD21F3} {\path{doi:10.1088/1674-1056/AD21F3}}.

\bibitem{Mao2023Physics-informedSolutions}
Z.~Mao, X.~Meng, {Physics-informed neural networks with residual/gradient-based adaptive sampling methods for solving partial differential equations with sharp solutions}, Applied Mathematics and Mechanics (English Edition) 44~(7) (2023) 1069--1084.
\newblock \href {https://doi.org/10.1007/S10483-023-2994-7} {\path{doi:10.1007/S10483-023-2994-7}}.

\bibitem{Liu2024AnPINNs}
Y.~Liu, L.~Chen, J.~Ding, Y.~Chen, {An Adaptive Sampling Method Based on Expected Improvement Function and Residual Gradient in PINNs}, IEEE Access (2024).
\newblock \href {https://doi.org/10.1109/ACCESS.2024.3422224} {\path{doi:10.1109/ACCESS.2024.3422224}}.

\bibitem{Tang2021DAS-PINNs:Equations}
K.~Tang, X.~Wan, C.~Yang, {DAS-PINNs: A deep adaptive sampling method for solving high-dimensional partial differential equations}, Journal of Computational Physics 476 (2021).
\newblock \href {https://doi.org/10.1016/j.jcp.2022.111868} {\path{doi:10.1016/j.jcp.2022.111868}}.

\bibitem{Tang2020DeepMapping}
K.~Tang, X.~Wan, Q.~Liao, {Deep density estimation via invertible block-triangular mapping}, Theoretical and Applied Mechanics Letters 10~(3) (2020) 143--148.
\newblock \href {https://doi.org/10.1016/J.TAML.2020.01.023} {\path{doi:10.1016/J.TAML.2020.01.023}}.

\bibitem{Wang2022IsNetwork}
C.~Wang, S.~Li, D.~He, L.~Wang, {Is L2 Physics-Informed Loss Always Suitable for Training Physics-Informed Neural Network?}, Neural Information Processing Systems (2022).
\newblock \href {https://doi.org/10.48550/ARXIV.2206.02016} {\path{doi:10.48550/ARXIV.2206.02016}}.

\bibitem{Rumelhart1986LearningErrors}
D.~E. Rumelhart, G.~E. Hinton, R.~J. Williams, {Learning representations by back-propagating errors}, Nature 323~(6088) (1986) 533--536.
\newblock \href {https://doi.org/10.1038/323533a0} {\path{doi:10.1038/323533a0}}.

\bibitem{Kingma2015Adam:Optimization}
D.~P. Kingma, J.~L. Ba, {Adam: A Method for Stochastic Optimization}, 3rd International Conference on Learning Representations, ICLR 2015 - Conference Track Proceedings (2015).

\bibitem{Hammersley1964MonteMethods}
J.~M. Hammersley, D.~C. Handscomb, {Monte Carlo Methods}, 1st Edition, Chapman and Hall, London, 1964.
\newblock \href {https://doi.org/10.1007/978-94-009-5819-7} {\path{doi:10.1007/978-94-009-5819-7}}.

\bibitem{Hinton2012LectureMagnitude}
G.~Hinton, T.~Tieleman, {Lecture 6.5 ‐ RMSprop: Divide the Gradient by a Running Average of Its Recent Magnitude} (2012).

\bibitem{Sutskever2013OnLearning}
I.~Sutskever, J.~Martens, G.~Dahl, G.~Hinton, {On the importance of initialization and momentum in deep learning}, in: Proceedings of the 30th International Conference on Machine Learning, 2013, pp. 1139--1147.

\bibitem{Kochenderfer2019AlgorithmsOptimization}
M.~J. Kochenderfer, T.~A. Wheeler, {Algorithms for optimization}, 1st Edition, The MIT Press, Cambridge, 2019.

\bibitem{Paszke2019PyTorch:Library}
A.~Paszke, S.~Gross, F.~Massa, A.~Lerer, J.~Bradbury, G.~Chanan, T.~Killeen, Z.~Lin, N.~Gimelshein, L.~Antiga, A.~Desmaison, A.~K{\"{o}}pf, E.~Yang, Z.~DeVito, M.~Raison, A.~Tejani, S.~Chilamkurthy, B.~Steiner, L.~Fang, J.~Bai, S.~Chintala, {PyTorch: An Imperative Style, High-Performance Deep Learning Library}, Neural Information Processing Systems (2019).

\bibitem{DelftHighPerformanceComputingCentreDHPC2024DelftBlue2}
{Delft High Performance Computing Centre (DHPC)}, \href{https://www.tudelft.nl/dhpc/ark:/44463/DelftBluePhase2}{{DelftBlue Supercomputer (Phase 2)}} (2024).
\newline\urlprefix\url{https://www.tudelft.nl/dhpc/ark:/44463/DelftBluePhase2}

\bibitem{Liu1989OnOptimization}
D.~C. Liu, J.~Nocedal, {On the limited memory BFGS method for large scale optimization}, Mathematical Programming 45~(1-3) (1989) 503--528.
\newblock \href {https://doi.org/10.1007/BF01589116} {\path{doi:10.1007/BF01589116}}.

\bibitem{Dauge1988EllipticDomains}
M.~Dauge, {Elliptic Boundary Value Problems on Corner Domains}, 1st Edition, Springer-Verlag, Berlin, 1988.
\newblock \href {https://doi.org/10.1007/BFb0086682} {\path{doi:10.1007/BFb0086682}}.

\end{thebibliography}


\end{document}